\documentclass[11pt,letterpaper]{amsart}
\usepackage{amsfonts, amssymb,amsmath, amsthm, amsxtra, latexsym, mathrsfs, amscd}
\usepackage{mathabx}
\usepackage{amsmath, amssymb, amsthm, times, float}
\usepackage{mathtools, amssymb, amsthm, times, float}
\usepackage{graphics}
\usepackage{MnSymbol}
\usepackage{tikz-cd}

\usepackage{amssymb}
\usepackage{hyperref}
\usepackage{backref}
\usepackage[latin1]{inputenc} %

\newtheorem{theo+}{Theorem}[section]
\newtheorem{prop+}[theo+]{Proposition}
\newtheorem{coro+}[theo+]{Corollary}
\newtheorem{lemm+} [theo+]{Lemma}
\newtheorem{lemm+def} [theo+]{Lemma and Definition}
\newtheorem{deep+}  [theo+]  {Deep Result}
\newtheorem{fact+}  [theo+]  {Fact}
\newtheorem{exam+}  [theo+]  {Example}
\newtheorem{rema+}  [theo+]  {Remark}
\newtheorem{defi+}  [theo+]  {Definition}
\newtheorem{xca+}[theo+]{Exercise}

\numberwithin{equation}{section}

\usepackage{color}

\newcommand\beq{\begin{equation}\label}
\newcommand\eeq{\end{equation}}

\renewcommand\({\Big(}
\renewcommand\){\Big)}
\renewcommand\[{\Big[}
\renewcommand\]{\Big]}
\newcommand\lgl{\langle}
\newcommand\rgl{\rangle}

\renewcommand\a[1]{{\acute{#1}}}

\newcommand\e[1]{{\ddot{#1}}}

\def\draft{\centerline{(Draft {\the \day}/{\the\month} \the \year.)}}
\bigskip

\def\refn#1.#2{\expandafter\def\csname#1\endcsname{[#2]}}
\def\refnr#1.{\csname#1\endcsname}

\def\fa{\mathfrak a}

\def\fe{\mathfrak e}
\def\ff{\mathfrak f}
\def\fg{\mathfrak g}

\def\fk{\mathfrak k}
\def\fl{\mathfrak l}
\def\fm{\mathfrak m}
\def\fn{\mathfrak n}

\def\fp{\mathfrak p}
\def\fq{\mathfrak q}

\def\ft{\mathfrak t}
\def\fu{\mathfrak u}

\def\fsl{\mathfrak {sl}}

\def\fso{\mathfrak{so}}
\def\fsp{\mathfrak{sp}}
\def\fspin{\mathfrak{spin}}
\def\fsu{\mathfrak{su}}

\def\a{\alpha}

\def\Claminv2{|C(\Lambda)|^{-2}}

\def\de{d\varepsilon}

\def\Aa2D{A^{\a,2}(D)}
\def\bAa2D{\overline{A^{\a,2}(D)}}
\def\Ab2D{A^{\beta,2}(D)}
\def\bAb2D{\overline{A^{\beta,2}(D)}}

\def\Norm#1_#2{\Vert#1\Vert_{#2}}

\def\phipl12{\phi_{p_{l_1}, p_{l_2}}}
\def\phip01{\phi_{p_{0}, p_{0}}}

\def\a{\alpha}

\def\Claminv2{|C(\Lambda)|^{-2}}

\def\sig{\sigma}

\def\ad{\operatorname{ad}}
\def\Ad{\operatorname{Ad}}
\def\Lie{\operatorname{Lie}}

\def\Cas{\text{Cas}}

\def\Ind{\operatorname{Ind}}

\def\exp{\operatorname{exp}}

\def\sgn{\operatorname{sgn}}

\def\tr{\operatorname{tr}}

\def\Ker{\operatorname{Ker}}

\def\rank{\operatorname{rank}}

\def\de{d\varepsilon}

\def\Aa2D{A^{\a,2}(D)}
\def\bAa2D{\overline{A^{\a,2}(D)}}
\def\Ab2D{A^{\beta,2}(D)}
\def\bAb2D{\overline{A^{\beta,2}(D)}}

\def\phipl12{\phi_{p_{l_1}, p_{l_2}}}
\def\phip01{\phi_{p_{0}, p_{0}}}

\def\alg/{algebra}

\def\Alg/{Algebra}

\def\alt/{alternative} 
\def\anal/{analytic}
\def\analfunc/{\anal/\ \func/}
\def\Ans/{\it Answer. \normal}
\def\ass/{associative}
\def\nass/{non-\ass/}
\def\autom/{automorphism}
\def\homom/{homomorphism}
\def\isom/{isomorphism}
\def\bdd/{bounded}
\def\Bdd/{Bounded}
\def\bddsymdom/{bounded \sym/ \dom/}
\def\Cartdom/{Cartan \dom/}
\def\bdry/{boundary}
\def\bsd/{\bdd/ \symdom/}
\def\bv/{boundary value}
\def\cf/{{\it cf}\.}
\def\Cf/{{\it Cf}\.}
\def\charr/{character}
\def\coeff/{coefficient}
\def\comm/{commutative}
\def\cpct/{compact}
\def\compl/{complex}
\def\comp/{complex}
\def\Comp/{Complex}
\def\conf/{conformal}
\def\conj/{conjugate}
\def\conn/{connect}
\def\cont/{continuous}
\def\conv/{converge} 
\def\convc/{convergence}
\def\convt/{convergent}
\def\convx/{convex}
\def\coord/{coordinate}
\def\lcoord/{local coordinate}
\def\Corr/{Corresponding}
\def\corr/{corresponding}
\def\corrd/{correspond}
\def\cov/{covariant}
\def\decomp/{decomposition}
\def\deco/{decompose}
\def\diff/{different} 
\def\Diff/{Different} 
\def\dimn/{dimension} 
\def\distr/{distribution} 
\def\div/{diverge} 
\def\dom/{domain}
\def\eg/{\hbox{\it e.g}\.}
\def\eigenf/{eigen\-\func/}
\def\eigensp/{eigen\-space}
\def\eigenv/{eigen\-value}
\def\eq/{equation}
\def\equa/{equation}
\def\de/{\diff/ial \equa/}
\def\do/{\diff/ial operator}
\def\ode/{ordinary \de/}
\def\pde/{partial \de/}
\def\pdo/{partial \diff/ial operator}
\def\psdo/{pseudo \diff/ial operator}
\def\fin/{finite}
\def\Ex/{\it Example.\ \normal}
\def\Exnr#1/{\it Example #1.\ \normal}
\def\foll/{follow}
\def\follg/{following}
\def\Follg/{Following}
\def\func/{function}
\def\Func/{Function}
\def\Fonc/{Fonc\-tion}
\def\fonc/{fonc\-tion}
\def\Funk/{Funk\-tion}
\def\funk/{Funk\-tion}
\def\gen/{general}
\def\har/{harmonic}
\def\Hint/{\it Hint. \normal}
\def\hist/{historic}
\def\histcl/{historical}
\def\hol/{holo\-morphic}
\def\homog/{ho\-mo\-ge\-ne\-ous}
\def\hyp/{hyper\-bolic}
\def\hyperg/{hyper\-geometric}
\def\ie/{\hbox{\it i.e.}}
\def\iff/{if and only if}
\def\ineq/{inequality}
\def\infra/{{\it inf\-ra}}
\def\ultra/{{\it ult\-ra}}
\def\Inpart/{In particular}
\def\inpart/{in particular}
\def\instof/{instead of}
\def\interps/{interpolation space}
\def\interp/{interpolation}
\def\Interp/{Interpolation}
\def\interpr/{Interpretation}
\def\Intr/{Introduction}
\def\intv/{interval}
\def\inv/{invariant}
\def\invc/{invariance}
\def\Iowords/{In other words}
\def\iowords/{in other words}
\def\ipr/{inner product}
\def\irred/{irreducible}
\def\lb/{line bundle}
\def\lin/{linear}
\def\lhs/{left hand side}
\def\rhs/{right hand side}
\def\loc/{local}
\def\math/{mathematic}
\def\mathcn/{\math/ian}
\def\manif/{manifold}
\def\meas/{measure}
\def\measl/{measurable}
\def\mero/{mero\-morphic}
\def\mon/{monomial}
\def\monog/{monogenic}
\def\mult/{multiple}
\def\multy/{multiply}
\def\multn/{multiplication}
\def\nas/{necessary and sufficient}
\def\nbd/{neighborhood}
\def\neg/{negative}
\def\nondeg/{nondegenerate}
\def\Oohand/{On the other hand}
\def\oohand/{on the other hand}
\def\Oonhand/{On the one hand}
\def\oonhand/{on the one hand}
\def\oper/{operator}
\def\orth/{ortho\-gonal}
\def\orthon/{ortho\-normal}
\def\otoh/{on the other hand}
\def\quat/{quaternion}
\def\pp/{\hbox{a. e.}}
\def\psh/{plurisubharmonic}
\def\pol/{polynomial}
\def\pot/{potential}
\def\pos/{positive}
\def\princ/{principle}
\def\prob/{probability}
\def\proj/{projective}
\def\projn/{projection}
\def\Proof/{\it Proof:\normal}
\def\Rem/{\it Remark\normal}
\def\Remnr#1/{\it Remark\ \normal #1. }
\def\rep/{representation}
\def\reps/{representations}
\def\meta/{metaplectic representation}
\def\repr/{reproducing}
\def\reprker/{reproducing kernel}
\def\resp/{respective} 
\def\resply/{respectively}
\def\restr/{restriction}
\def\sa/{self-adjoint}
\def\st/{such that}

\def\sol/{solution}
\def\ru/{space}
\def\sph/{spherical}
\def\ssp/{sub\ru/}
\def\sym/{symmetric}
\def\Sym/{Symmetric}
\def\symb/{symbol}
\def\symbc/{symbolic}
\def\symdom/{\sym/ domain}
\def\symp/{symplectic}
\def\Theor#1/{\fet Theorem #1.\ \normal}
\def\Lem#1/{\fet Lemma #1.\ \normal}
\def\Lemma/{\fet Lemma.\ \normal}
\def\topl/{topology}
\def\topll/{topological}
\def\transf/{transform}
\def\transl/{translation}
\def\transfn/{transformation}
\def\transv/{transvectant}
\def\trig/{trigonometric}
\def\tril/{trilinear}
\def\trilf/{trilinear form}
\def\uhp/{upper halfplane}
\def\uhs/{upper halfspace}
\def\vb/{vector bundle}
\def\vf/{vector field}
\def\vsp/{vector space}
\def\wrt/{with respect to}
\def\Wlog/{Without loss of generality}
\def\a{\alpha}

\def\e{\varepsilon}

\def\sig{\sigma}

\def\Ab/{Abel}
\def\Ban/{Banach}
\def\Bansp/{\Ban/ space}
\def\Belt/{Bel\-tra\-mi}
\def\Berg/{Berg\-man}
\def\Bern/{Ber\-nou\-lli}
\def\Berz/{Berezin}
\def\Bess/{Bessel}
\def\Cart/{Car\-tan}
\def\Cay/{Cay\-ley}
\def\CG/{Clebsch-Gordan}
\def\Cl/{Clifford}
\def\CR/{Cauchy-Rie\-mann}
\def\Dir/{Dirichlet}
\def\Eucl/{Euclide}
\def\Eucln/{Euclidean}
\def\F/{Fourier}
\def\Hank/{Hankel}
\def\Hankf/{\Hank/ form}
\def\Herm/{Hermite}
\def\Hilb/{Hilbert}
\def\Hilbs/{Hilbert space}
\def\Hilbsp/{Hilbert space}
\def\HS/{Hilbert-Schmidt}
\def\Lag/{La\-grange}
\def\Lap/{La\-place}
\def\LapBelt/{\Lap/-\Belt/}
\def\Leb/{Lebesgue}
\def\Marc/{Mar\-cin\-kie\-wicz}
\def\Moeb/{Moebius}
\def\Moebt/{Moebius transformation}
\def\Moebtransfn/{Moebius transformation}
\def\Pla/{Plan\-che\-rel}
\def\Poin/{Poin\-car\'e}
\def\Riem/{Rie\-mann}
\def\Riemn/{\Riem/ian}
\def\psRiemn/{pseudo-\Riem/ian}
\def\Riems/{Rie\-mann surface}
\def\Schroe/{Schr\"odinger}
\def\Weier/{Weier\-strass}
\def\anal/{analytic}
\def\bsd/{bounded symmetric domain  }
\def\bdd/{bounded}
\def\calc/{calculation}\def\conj{conjugate}
\def\calci/{calculating}\def\eg{e.g.}
\def\conj/{conjugate}
\def\deco/{decomposition}
\def\eg/{e.g.}
\def\fct/{function}
\def\gp/{group}
\def\hw/{highest weight}
\def\hwv/{highest weight vector}
\def\hwvs/{highest weight vectors}
\def\lw/{lowest weight}
\def\lwv/{lowest weight vector}
\def\lwvs/{lowest weight vectors}
\def\hds/{holomorphic discrete series}
\def\iff/{if and only if}
\def\inv/{invariant}
\def\irrde/{irreducible decomposition}
\def\meas/{measure}
\def\transf/{transform}
\def\rep/{representation}
\def\resp/{respectively}
\def\inters/{intertwines}
\def\interg/{intertwining}
\def\meta/{metaplectic representation}
\def\qu/{quaternion}
\def\rep/{representation}
\def\symdom/{ symmetric domain}
\def\st/{such that}
\def\shd/{subhead}
\def\transf/{transform}
\def\wrt/{with respect to}

\def\Norm#1#2#3{\Vert#1\Vert^{#3}_{{#2}}}

\def\tr{\operatorname{tr}}

\begin{document}

\title[
Principal series 
induced from  Heisenberg parabolic subgroups
]
{Principal series of quaternionic
  and real
split exceptional Lie groups
 induced from  Heisenberg parabolic subgroups
}
\author{Genkai Zhang}
\address{Mathematical Sciences, Chalmers University of Technology and
Mathematical Sciences, G\"oteborg University, SE-412 96 G\"oteborg, Sweden}
\email{genkai@chalmers.se}
\thanks{ Research supported partially
 by the Swedish
Research Council (VR, Grants 2018-03402, 2022-02861).
}
\begin{abstract}
  Let $G$ be an irreducible quaternionic Lie groups rank $4$
  or an exceptional real split
  Lie group $E_{6(6)}$,
$E_{7(7)}$, or $E_{8(8)}$
  We study the
principal series
representation $\pi_\nu=\Ind_P^G(1\otimes e^\nu\otimes 1)$
of $G$
induced   from
  the Heisenberg parabolic subgroup
  $P=MAN$ realized on $L^2(K/L)$, $L=K\cap M$.
  We find the $K$-types in
  the induced representation
via a double cover $K/L_0\to K/L$
  and a circle bundle
  $K/L_0\to K/L_1$ over
  a compact Hermitian symmetric space $K/L_1$.
We compute
the Lie algebra $\fg$-action of
$G$ on the representation space. We
find the complementary series,  reducible points, and unitary
  subrepresentations   in this family of representations.
\end{abstract}

\maketitle
\tableofcontents

\baselineskip 1.40pc

\section{Introduction}

The present paper is  a continuation
of our  earlier paper
  \cite{Z22} where we studied
  the   representations
  of Hermitian Lie groups $G$
 induced
 from Heisenberg parabolic subgroups.
  Here we consider  the case of quaternionic 
  Lie groups $G$ of real rank $4$
  and split exceptional Lie groups of types $E_6, E_7, E_8$.  We shall find
  the complementary series, the reducible points and subrepresentations.
  
\subsection{The induced representations}  
Let $G$ be an irreducible quaternionic
Lie group or the split exceptional Lie group $E_{6(6)}, 
E_{7(7)}, E_{8(8)}$; the maximal
compact subgroup  $K$  of $G$ is 
locally
isomorphic  $SU(2)\times M_c$, or 
$Sp(4), SU(8), Spin(16)$
in the respective cases.
Let  $\fg=\fk+\fp$ be the Cartan decomposition
of $\fg=\Lie(G)$,   $\fk=\Lie(K)$.
 The rank of $G/K$ is  $r=1, 2, 3, 4, 6, 7, 8$
and we shall treat only the case $r=4$
and the exceptional cases
with $r= 6, 7, 8$. 
The spherical principal series
for the rank   one quaternionic group  $G=Sp(n, 1)$
 has  been well-studied \cite{JW}.

We  choose
  a  Cartan subalgebra in $\fk^{\mathbb C}$
and extend it to a Cartan subalgebra in   $\fg^{\mathbb C}$;
all symmetric pairs
$(\fg, \fk)$ are of  of equal rank,
$\rank \fk^{\mathbb C}=\rank \fg^{\mathbb C}
$ except for
$(\fg, \fk)=(\fe_{6(6)}, \fsp(4)$,
and  Cartan subalgebra is aso a Cartan subalgebra for
 $\fg^{\mathbb C}$.
 Let $\beta_1$ and $\gamma$ be the highest compact
 and non-compact roots, respectively. 
Then $\beta_1$   defines a $\fsu(2)$-ideal in $\fk=\fsu(2) +\fm_c$
 when $(\fg, \fk)$ is quaternonic. The
   real 
   part $E=E_\gamma + E_{-\gamma}\in \fp$ of
the   root vector $E_\gamma$ of 
   $\gamma$  defines
   a $\fsl(2, \mathbb R)$-subalgebra
  and a  Heisenberg parabolic subalgebra
  $\fm + \fa +\fn
  $ in $\fg$, $\fa=\mathbb R E$.
  Let $P=MAN$ be the corresponding Lie subgroup of $G$.
  We shall study the induced representations
  $\pi_\nu=\Ind_{P}^G(1\otimes e^\nu\otimes 1)$
  and $  \pi_{\nu, -} =\Ind_{P}^G(\sgn\otimes e^\nu\otimes 1)$
  of $G$ and find their  complementary series and study their composition series.

  The representation   $\pi_\nu$ of $G$
  can be realized in the so-called compact picture 
  on the space $L^2(K/L)$,  $L=K\cap M$, 
  with the action of $G$ restricted to $K$ being  the left regular action of $K$.
  The structure of the
  representation
  of $\pi_\nu$
  is determined by
  the action of $\fp$ on $L^2(K/L)$. 
  Generally it is rather challenging to 
  find the action, and the most studied
  case is when $K/L$ is a symmetric space with $L^2(K/L)$
  being of multiplicity free.
  In the present case $L^2(K/L)$ 
  under $K$ is 
 of arbitrarily high multiplicities   and we shall  compute the action by developing further
  our method in \cite{Z22}
 when $G$ is  Hermitian.

\subsection{The main results and  methods}  

We realize $G$ as the adjoint group
of $\fg$. There is a double cover
$K/L_0\to K/L$ of $K/L$
and a circle $U(1)$-bundle  $K/L_0\to K/L_1$ 
over a compact Hermitian 
  symmetric space $K/L_1$
  or rank four with $K/L_1=  \mathbb P^1  \times M_c/L_c$
  being reducible for quaternionic $G$, $L_c=M_c\cap L_1$,
  with $\mathbb P^1=SU(2)/U(1)$ the projective space
  and $M_c/L_c$ of rank $3$,
  and $K/L_1$ being irreducible
  for the exeptional pairs $G=E_{6(6)},
  E_{7(7)},E_{8(8)}  $.
The space $L^2(K/L)$ consists of 
 even functions in $L^2(K/L_0)$
which can  expanded
using the Fourier series
  along the fiber  $U(1)$ with coefficients being  $L^2$-sections of line bundles over $K/L_1$.
 The $L^2$-space of sections 
 can further be decomposed 
  using the Cartan-Helgason-Schlichtkrull
  theorem \cite{Sch}. We  compute
  then  the Lie algebra action of $\fa
=\mathbb R E$ on the
  $L$-invariant element in 
  each irreducible subspace in $L^2(K/L)$.
  This is done by developing
recurrence formulas for products
and
differentiations
  of spherical polynomials;
the product formulas
is found by using Harish-Chandra
$c$-functions \cite{Vretare, Z-tams, Z22}  and  the differentiations
by using the techniques in \cite{B-O-O}.

Our first main result is Theorem    \ref{pi-E-W}
where we compute the Lie algebra action.
As a consequence we get the range
of complementary series,  Theorem \ref{compl-ser},
and the
reducible points, Theorem \ref{redu-pt}.
We prove in Proposition \ref{min-as-sub}
that the minimal representation
is a proper subrepresentation of the kernel
of the system conformally invariant
differential operators introduced
in \cite{B-K-Z, GW}.

  There are some  geometric differences
  between the representations $\pi_\nu$
  for our present cases (with 
  $\fg$ being quaternionic or $\fe_{6(6)},
  \fe_{7(7)},   \fe_{8(8)}$)
  and for Hermitian Lie groups treated in \cite{Z22}. When $G$ is
a  Hermitian Lie group the subgroups $M$ and  $L=M\cap K$ are
  connected. This is due to the fact that $K$
  has one-dimensional center and $\fp^{\mathbb C}=\fp^+ +\fp^-$
  is reducible, and that any group element in $K$ fixing the element
  $E=E_\gamma +E_{\gamma}$ also fixes
  $  E_{\pm \gamma} $.
    For quaternionic Lie groups
or   $E_{6(6)}/E_{7(7)}/E_{8(8)}$
  we have two
  families of induced representations,
  one induced from the trivial representation
of $M$  and one from the sign representation of $M$.
The construction of $K$-finite elements in $L^2(K/L)$
involves some averaging process of
the elements
obtained from
the Cartan-Helgason-Schlichtkrull theorem
for $K/L_1$.

\subsection{Brief account of earlier results and future questions}  
The study of induced
representations for quaternionic Lie groups 
from Heisenberg parabolic subgroups
is also motivated from the theory of
minimal representations.
Indeed  it has been proved by Gross and Wallach \cite[Sections 13-14]{GW}
  that the minimal quaternionic
  representations for $E_{6(2)}, E_{7(-5)}, E_{8(-24)}$
and $ SO(4, 4)$
  can be realized as subrepresentations of
  the  principal series;
 they appear also in the analytic continuation
 of the quaternionic discrete series.
 As an application of our $K$-type formula
 we prove that the minimal
 representation is a proper sub-representation
 in the kernel of some system
 of conformally invariant differential
 operators for
 $ E_{6(2)}, E_{7(-5)}, E_{8(-24)}$
 and $SO(4, d), d\ge 4 $,  $d$ being even.
 There leaves still the question
 of understanding the full composition series, and
due to the higher multiplicities it looks
 still elusive.

There has been much
study on complementary series and composition series
for induced representations and we give
here a very brief account of the most
relevant known results.
The Heisenberg parabolically induced representations  
for $G=SU(p, q)$
    and similar representations for $SO(p, q)$
    are studied in  \cite{HT}
    and in  \cite{FWZ, Z22}
for general Hermitian Lie groups.
 Induced representations $\Ind_P^G(1\otimes e^\nu\otimes 1)$
 from maximal parabolic subgroups of
 a large class of Lie groups $G$
 with  $L^2(K/P\cap K)$ being multiplicity
 free have been studied extensively;
 see e.g. \cite{H,  Lee, Lee-Loke, Lee-Zhu, OZ-duke, Sahi, Z95}
 and references therein.
 
For Hermitian Lie groups 
the  induced representations
are also studied using Fourier analysis on Heisenberg
groups $\bar N=N_{-1}\times N_{-2}=N_{-1}\times \mathbb R$
in \cite{FWZ}.
This involves the decomposition of
the metaplectic representation
of $
Sp(n, \mathbb R)
$, $n=\frac 12 \dim N_{-1}$, 
under subgroup $M$.
This decomposition is discrete
as the Hermitian symmetric spaces of $M$
and $Sp(n, \mathbb R)$
have the same complex structure.
It might
be interesting to study
the induced representation
for quaternionic Lie groups
using also the non-compact realization.
However in this case the decomposition
above is not discrete
\cite{Z-22-1} and
more analysis problems are involved.
I was also informed that
some results here maybe obtained
from the Atlas software (http://www.liegroups.org/software/).
The systematic approach here using
geometry and analysis of the circle bundles
might be of interests in its own right.

In a forthcoming work \cite{FHWZ-g2}
we shall study the maximal parabolically
induced representations for the group $G_2$;
there are two families of maximal
parabolic subgroups one of which is Heisenberg.
This case is somewhat different from the 
rank four cases. There
are also Heisenberg parabolic subgroups
for the exceptional split real groups
and for their complexification.
In \cite{YWZ}
we have found a generalization
for the Cartan-Helgason-Schlichtkrull theorem
for quaternionic twistor spaces
and this will enable us to study
a larger class of induced representations.

\subsection{Organization}  

  The paper is organized as follows. In Section
  \ref{prel-sect}
  we recall
  some basic facts about
  the quaternionic Lie algebras
  and exceptional real forms $e_{6(6)}/e_{7(7)}/e_{8(8)}$
  and strongly orthogonal compact and non-compact roots.
  In Section \ref{h-section}
  we  construct
 the  Heisenberg  parabolic group
  $MAN$ and study the geometry of $K/L$
  via circle bundles over
  compact Hermitian symmetric spaces.
  The induced representations of $G$
are introduced and
their  $K$-types are found
Section \ref{deco-sect}
via its realization in $L^2(K/L)$
and the coverings.
 Section \ref{main-sect}
  is devoted to the Lie algebra action on
  $L^2(K/L)$.
  The  complementary series are studied in Section \ref{comp-sect}.
  The relation between minimal
  representations of quaternionic Lie groups
  and kernels some of conformal invariant
differential operators  is studied in   Section \ref{quat-ker}. In the appendix we prove some  product formulas for spherical polynomials for line bundles
over compact Hermitian symmetric spaces.

  \subsection*{Acknowledgements} I would like to
  thank Dan Barbasch, Jan Frahm
    and Clemens Weiske
    for several inspiring
    discussions.

  \subsection*{Notation}

\begin{enumerate}
  \item $(\fg, \fk)$: a  quaternionic
  symmetric pair of rank $4$ or exceptional pair
  $            (\fe_{6(6)}, \fsp(4))$,
$              (\fe_{7(7)}, \fsu(8))$,
$                            (\fe_{8(8)}, \fso(16))$
                            or rank $ 6, 7, 8$.

  \item     $\fg=\fk+\fp
    $:
  Cartan decomposition of $\fg=\Lie(G)$
  with Cartan involution $\theta$;
    $\fk =\fsu(2)+
  \fm_c$; $\ad X (Y)=[X, Y]$: adjoint representation in $\fg$;
$\Ad g (h)=g hg^{-1}$, $gX=\Ad g (X), X\in \fg$: adjoint actions
in $G$ on $G$ and on $\fg$.

\item     $\fg^{\mathbb C}=\fk^{\mathbb C}+\fp^{\mathbb C}$:
  complexification of $\fg$, with the corresponding
  Cartan involution $\theta$;
  $X\to \overline{X}$: conjugation in $\fg^{\mathbb C}$
  with respect to the real form $\fg$;
$B(X, Y)=(X, Y)$: the Killing form on $\fg^{\mathbb C}$ normalized
with the highest root $\beta_1$ having dual square norm
$(\beta_1, \beta_1)=2$.

\item $ \ft^{\mathbb C}=\mathbb C {\beta_1}^{\vee}  + (\ft\cap\fm_c)^{\mathbb C}
  \subset \fsl(2)_{\beta_1}+ \fm_c^{\mathbb C}$:
  Cartan subalgebras
  of $\fk^{\mathbb C}$ and
${\fg}^{\mathbb C}$.

\item   $\gamma$: highest non-compact root
  of $\fg^{\mathbb C}$.

\item   $E_{-\alpha} =-\overline{\theta(E_{\alpha})}$: 
  choice of 
  $E_{-\alpha}$ with a given $E_{\alpha}$
for a positive root $\alpha$.

\item $\{E_{\beta_1}, E_{-\beta_1}, {\beta_1}^{\vee}\}$:
  $\fsl(2)$-triples in $\fg^{\mathbb c}$;
   $\{E_{\beta_1}+ \overline{E_{\beta_1}}, i(E_{\beta_1}- \overline{E_{\beta_1}}),  i{\beta}^{\vee}\}$:
  $\fsu(2)$-triples in $\fg$.

\item  
  $\gamma_1, \gamma_2, \gamma_3, \gamma_4$:
  strongly orthogonal roots in
  $\fp^{\mathbb C}$,
  with co-roots
  $H_1, H_2, H_3, H_4$,
  $\gamma_i(H_i) =2, \gamma_i(H_j)=0,  1\le i, j\le r$.
  
\item
$
  \beta_1=\gamma + \gamma_1$:
relation between highest compact and non-compact roots and lowest non-compact root;
  $\beta_j=
\gamma -\gamma_j =\gamma -\gamma_1-\gamma_j, \,
2\le j\le 4
$: strongly orthogonal compact roots.

  
\item $\mu=(\mu_1, \mu_2, \mu_3, \mu_4)=
 \frac 12( \mu_1 \beta_1
  + \mu_2
\beta_2 +\mu_3\beta_3+  \mu_4 \beta_4)$: dominate  weights
  for the symmetric pair
  $(\fk^\ast, \fl_1)$, spherical
  with respect to one-dimensional
  representation of $\fl_1$; $\rho=\rho_{\fk^{\mathbb C}}$: half-sum of positives roots in 
$\fk^{\mathbb C}$.

\item       $\fa=\mathbb RE \subset \fp$,
  $E=E_\gamma +E_{-\gamma}$:
 one-dimensional  subalgebra defining the Heisenberg parabolic subalgebra 
  $\fm +\fa +\fn$. 

\item  $\fl =\fk\cap \fm$.

\item $\fn=\fn_1+\fn_2$,  
$\fn_2= \mathbb R(F- iH) 
=\{X\in \fg; [E, X]=2X\}$, $
\fn_1= \{X\in \fg; [E, X]=X\}.$  

\item $\fk=\fk_2+\fk_1+\fl$: orthogonal
  decomposition of $\fk$ obtained from 
  $\fg= (\fn_2 +\fn_{-2}) +
  (\fn_1 +\fn_{-1}) + (\fm +\fa)$.

\item $P=MAN$:
  parabolic subgroup; $L=K\cap M$.

\item $L_0$: connected component of $L$;
  $L_1\subset K$: stabilizer of the line $[\mathbb CE_\gamma]\in
  \mathbb P(\fp^{\mathbb C})$
  in the projective space.

  \item $(\fk, \fl_1)$: compact Hermitian symmetric pair
  of $\fk$ of rank $4$.

\item  
$\chi^{p}$:
  character on $L_1\to L_1/L_0=U(1)\to U(1)$, $p\in \mathbb Z$; $L^2(K/L_1; \chi^p)$: $L^2$-space
  of sections of the line bundle on $K/L_1$ defined by $\chi^p$.

\item $w_0\in K$:  element in $L=M\cap K$, $w_0^2= 1$,
  such that $L=L_0\rtimes \{1, w_0\}$
  and $\Ad(w_0)|_{\fn_{\pm 2}} = -1$.

\item $\phi_{\mu; l,p}(k)$: spherical
  polynomial for $K/L_1$ of type
  $(l, p)$;
  $\phi_{\mu; l}(k)=\phi_{\mu; l,l}(k)$: spherical
  polynomial for $K/L_1$ of type
  $(l, l)$.

\item $\Cas=\Cas_0 +\Cas_1+\Cas_2$: Casimir element
  of $\fk$ written as
  a sum of Casimir elements under
  the decomposition  $\fk=\fk_2+\fk_1+\fl$;
  $\Cas(\mu) 
  = -(\mu + 2\rho_{\fk^*}, \mu)$:
  eigenvalue of the Casimir operator of $\fk^{\mathbb C}$
acting on the representation with highest weight $\mu$.
\end{enumerate}

\section{Quaternionic Lie algebras
  and split real Lie algebras $e_{6(6)},
  e_{7(7)}
  e_{8(8)}$
}\label{prel-sect}

\subsection{Simple Lie algebras  $\fg$
  and their $\fsl_2$-subalgebras}
\label{sect-1.1}
Let $\fg$ be a simple real {\it non Hermitian} Lie algebra
and $\fg=\fk +\fp $ its Cartan decomposition
with respect to the Cartan involution $\theta$.
Let $\fg^{\mathbb C}$ 
be  its complexification,
$\fg^u=\fk +i \fp $ the compact real form
of $\fg^{\mathbb C}$,   $G^{\mathbb C}=\Ad(\fg^{\mathbb C})$  the
 {\it adjoint group} of $\fg^{\mathbb C}$,
and $G=\Ad(\fg)$
the corresponding adjoint subgroup
of  $\fg$. Let $K$ be the maximal compact
subgroup of $G$. In particular $K$
and $G$ are all connected. To save notation
the adjoint action of $G$ on 
$\fg^{\mathbb C}$ will be just written
as $gX:=\Ad(g) X$, $g\in G, X\in \fg^{\mathbb C}$
when no-confusion would arise.
This realization of  $G$
and $K$ has 
some  advantage
of explicit description of certain
homogeneous spaces  of $K$ as
submanifolds of the linear space
$\fp^{\mathbb C}$.

We fix first the standard  $\fsl(2, \mathbb C)$-triple,
  $$
 H=\begin{pmatrix} 1&0 \\
  0 & -1 
\end{pmatrix},\, 
E_{+} =
\begin{pmatrix} 0&1 \\
  0 & 0 
\end{pmatrix}, \, 
E_{-} =
\begin{pmatrix} 0&0 \\
  1 & 0 
\end{pmatrix}, 
  $$
  with the Killing form $B(\cdot, \cdot)$
  on $\fg^{\mathbb C}$ normalized by $B(E, E)=\tr E E=2$.
As a convention the  compact form  $\fsu(2)$
and non-compact real form and their Cartan decompositions 
will be chosen as 
\begin{equation}
\label{eq:sl2-su2}
\fsu(2) 
 =\mathbb RiH +
( \mathbb Ri E+
\mathbb R iF)\subset \fg^u, \,
\fsu(1, 1)=\mathbb RiH +
(\mathbb R E + \mathbb R F)  \subset \fg,
\end{equation}
where
$$E=E_+ + E_-, \quad  F= i(E_+-E_-), \quad E_+=\frac 12(E-iF),
\quad E_+=\frac 12(E+iF).$$

We choose
 $iH\in \fsu(2) $
or $E\in \fsl(2, \mathbb R)$
as Cartan element of
the respective subalgebras
with $iH$ having imaginary roots.
  The Cayley transform 
  $$
\Ad\(\exp 
(-\frac{\pi i}{4}
(E_+- E_{-}))
\)=
\Ad\begin{pmatrix}
  \cos\frac{\pi}4 &-\sin\frac{\pi}4 
  \\
\sin\frac{\pi}4   &  \cos\frac{\pi}4 
\end{pmatrix}:  \fg^{\mathbb C}\to  \fg^{\mathbb C}, \,
H\mapsto   E_+ + E_-
$$
exchanges two simple elements $H$ and $E_++E_-$.
For later references we recall also 
the elementary fact  that
\begin{equation}
\label{-1-in-su2}
\exp(\pi i H)=\begin{pmatrix} -1&0 \\
  0 & -1 
\end{pmatrix}=:-1. 
\end{equation}

Let  $\fsl(2, \mathbb C)
\subset \fg^{\mathbb C}$
be a $\theta$-stable subalgebra. If it is determined by
a compact root $\alpha$ (see below)
with co-root ${\alpha}^{\vee}
\in
\fsl(2, \mathbb C)$,  ${\alpha}(\alpha^{\vee})=2$,  we shall denote it by
$\fsl(2, \mathbb C)=\fsl(2, \mathbb C)_{\alpha}\subset \fg^{\mathbb C}$
and the compact subalgebra $\fsu(2)_{\alpha}\subset \fg$
as in (\ref{eq:sl2-su2}).
We  write sometimes
the element 
$-1\in SU(2)=SU(2)_\alpha$
in (\ref{-1-in-su2}) as
$(-1)_\alpha$ to indicate the dependence on $\alpha$, as well as
$U(1)_{{\alpha}}:=\exp(\mathbb Ri{\alpha}^{\vee})
\subset K\subset G$.
Finite dimensional irreducible representations
of
$\fsu(2)_{\alpha}$
will be parametrized as $\frac m2 \alpha$
identified with integers $m\ge 0$. The Casimir
element $\Cas=\Cas_{\alpha}$ of
$\fsu(2)_{\alpha}$ is $\frac 12 H^2 + E_+ E_- + E_- E_+$
and has  the eigenvalue
$-\frac 12 (m+2)m$ on the representation space. If  $\fsl(2, \mathbb C)
\subset \fg^{\mathbb C}$ is determined by a real root $\alpha$ the corresponding real form in $\fg$
will be written as 
$\fsl(2, \mathbb \mathbb R)_\alpha$.

\subsection{Quaternionic Lie algebras $\fg$}
Let $(\fg, \fk)
=(\fk+\fp, \fk)$ be
 an irreducible quaternionic symmetric pair of real rank $r$.
  Then  $r=1, 2, 3, 4$. If $r=1, 2$ then 
$\fg=\fsp(1, n), \fsu(2, n), \fg_2$; if $r=3$ then 
$\fg= \fso(4, 3)$. 
We shall only consider
the case $r=4$;  a list of all such $\fg$ is given in the Table \ref{tab:1}; see also \cite{GW} for a complete list.
By definition   $\fk =\fsu(2) +\fm_c$ is a sum of two
  ideals.
  The space $\fp^{\mathbb C}=\mathbb C^2\otimes V$,
  as representation of
  $\fk =\fsu(2) +\fm_c$ with  $V$
  an irreducible representation of $\fm_c$.

  Choose  $\ft=\ft_{\fsu(2)}
  +\ft_c\subset \fk=\fsu(2)+\fm_c\subset \fg$  
   Cartan subalgebra of $\fk$ and $\fg$.    We let $R=R(\fg^{\mathbb C}, 
   \ft^{\mathbb C}) 
   $
   be the set of   roots 
   and 
 $\fg^{\mathbb C}_\alpha$
 the corresponding root spaces, $\alpha\in R$, 
 with $R^+ $ the positive roots with respect to a fixed 
 ordering.
 Let 
  $\Delta=\Delta(\fg^{\mathbb C}, 
  \ft^{\mathbb C})$
  be the simple roots. 
    Then  $R^+=R_c^+\cup R_n^+$, 
  the compact and  non-compact 
positive roots, respectively, and 
  $\Delta$ contains a unique 
  non-compact root $\gamma_1$. 
 Recall that the  co-root ${\alpha}^{\vee}$ of a root
 $\alpha$ is defined by  $\beta(  {\alpha}^\vee)=\frac {2(\beta, \alpha)}{(\alpha, \alpha)}$
for all $\beta\in \Delta$. 
 The ordering is chosen so that
the highest root   $\beta_1$  is a compact root and
the
corresponding $\fsu(2)_{\beta_1}$-subalgebra 
is the ideal $\fsu(2)$ in $\fk=\fsu(2)+\fm_c$ above, so that
the Cartan subalgebra $\ft^{\mathbb C}_{\fsu(2)}$
is generated by the co-root  ${\beta_1^\vee}$,
$\ft^{\mathbb C}_{\fsu(2)}=\mathbb C {\beta_1^\vee}
=\mathbb C H_{\beta_1}$ as in 
Section \ref{sect-1.1}.

  Thus  the compact positive roots are 
 $ R_c^+ = \{\beta_1\} \cup 
 R^+(\fm_c^{\mathbb C}, \ft_c^{\mathbb C} ) 
 $.

 Let $\gamma$  be the highest non-compact root, i.e.
 the highest weight of
 $\fp^{\mathbb C}=\mathbb C^2\otimes
 V$
 as representation of $\fk^{\mathbb C}=\fsl(2, \mathbb C)_{\beta_1}
 +\fm_c^{\mathbb C}$; a precise
 formula will be given in Lemma \ref{lem1-1}.
 We shall fix from now on the  Killing form
$B(\cdot, \cdot):=(\cdot, \cdot)$ on $\fg^{\mathbb C}$
so that 
the highest root  $\gamma$
has square norm $2$, $
(\gamma, \gamma)=2$.
This normalization implies
also that
$(H, H)=
(\gamma, \gamma)=
2 $
for $H:=\gamma^{\vee}$.
By  orthogonal decomposition 
of subspaces of  $\fg^{\mathbb C}$ we always 
mean with respect to the Killing form 
$B(\cdot, \cdot)$.
 An Hermitian 
inner product on a vector space $V$ will 
be written as $\langle \cdot, \cdot\rangle$
and we fix the corresponding Hermitian inner product
on $\fp^{\mathbb C}$, and the positive definite
Hermitian inner product on $\fk^{\mathbb C}$
by extending the positive
definite form $-B(\cdot, \cdot)$ on $\fk$.
We fix generators
for $\fsl(2, \mathbb C)_{\gamma}
=\mathbb C H+\mathbb C E_{\gamma} +
\mathbb R E_{-\gamma}$ with the Killing form
$$
B(E_\gamma, E_{-\gamma})
=\langle E_\gamma, E_{\gamma}\rangle
=1.
$$
The real form $\fsl(2, \mathbb R)_{\gamma}\subset \fg$ is 
$
\fsl(2, \mathbb R)_{\gamma}=\mathbb R E+
\mathbb R F+ \mathbb R iH$,
  \begin{equation}
    \label{X0-Y0}
    E= E_{\gamma} + E_{-\gamma}
    =E_{\gamma} +\overline{E_{\gamma}}, 
       F = i(E_{\gamma} - E_{-\gamma}),  H={\gamma^\vee}. 
     \end{equation}
     and $B(E, E)=B(F, F)=2$,
     with the commutation relation as in (\ref{eq:sl2-su2}).

     We recall the following fact
     on the realization of $K$.

\begin{lemm+}   (\cite{GW}.) 
  Under the realization of $K$
  via the adjoint action
  on $\fp^{\mathbb C}=\mathbb C^2\otimes
 V$ it is  $K=SU(2)_{\beta_1} \times M_c/
((-1)_{\beta_1}\times\e)$, where
$M_c$ is a compact  subgroup
of $GL(V, \mathbb C)$
with central element $\e$
and Lie algebra $\fm_c$.
\end{lemm+}
 
The group $K$ will always be realized
as linear maps  on $\fp$ and $\fp^{\mathbb C}$, and
finite dimensional
representations of $K$
will  be realized as representation
of $SU(2) \times M_c$
with trivial action by $(-1)\times \e$.
 The space $SU(2)/U(1)
= SU(2)_\beta/U(1)_\beta$
  is the projective space $\mathbb P^1$, 
  $U(1)=\exp(i\mathbb R{\beta^\vee})$.

\subsection{Strongly orthogonal roots
  $\beta_1, \cdots, \beta_4$ and 
  $\gamma_1, \cdots, \gamma_4$
 for  quaternionic Lie algebra $\fg$}
\label{strong-ort}
Starting with the simple non-compact
  root $\gamma_1$
  we can find
   a system $\gamma_1<\cdots <\gamma_4$
 of strongly orthogonal 
 non-compact roots. Similarly
 starting with the highest
 compact roots $\beta_1$
 there is a system
 $\beta_1>\cdots >\beta_4 $
 of strongly orthogonal compact roots.
 We shall need some relations
 between these two systems.

 When $\fg=\fso(4, d), d\ge 4$
 we have $\fk=\fsu(2)+
 \fsu(2)+ \fso(d)$ with $\fm_c=
 \fsu(2)+ \fso(d)$ is reducible.
 This case is
 slightly different from the other
 cases  so we work out the details here.

 The root system of $\fg^{\mathbb C}=\fso(4+d, \mathbb C)$
 is of type $B$ or $D$.
 Fix  $\{\e_i, j=1, \cdots, [\frac d 2]\}$ an orthonormal basis of    $\mathbb R^{[\frac d 2]}$.
  In the usual notation the root system of 
  $\fg^{\mathbb C}$ is 
  $$
  R_n^+=\{\e_1, \e_2\}\cup 
  \{\e_1\pm \e_j, \e_2\pm \e_j, 3\le j\le \frac{
  d-1}2 
  \}, \, 
  $$
  $$
  R_c^+=\{\e_1\pm \e_2\}\cup 
  \{\e_i\pm \e_j, 1\le i< j\le \frac{
  d-1}2 
\}
\cup 
  \{\e_j, 2< j\le \frac{
  d-1}2 
\}
,\, 
  $$
  if $d$ is odd, and
$$
  R_n^+=  \{\e_1\pm \e_j, \e_2\pm \e_j, 1\le j\le \frac{
  d}2 
  \}, \, 
  R_c^+=\{\e_1\pm \e_2\}\cup 
  \{\e_i\pm \e_j 1\le i< j\le \frac{
  d}2 
  \},
  $$
  if $d$ is even. The highest
  non-compact root is $\gamma=\e_1+\e_3$.

The strongly orthogonal 
strongly  compact roots 
$\{
\beta_1, \beta_2, \beta_3, \beta_4\}$
and  non-compact roots $\{\gamma_1, \gamma_2,  \gamma_3, \gamma_4\}$
are 
$$
\beta_1=\e_1+\e_2>\beta_2=\e_1-\e_2 >\beta_3=
\e_3+\e_4 >\beta_4=\e_3-\e_4,
    $$
    $$
  \gamma=\e_1+\e_3>
    \gamma_4=\e_1+\e_4 >
    \gamma_3= \e_1-\e_4 >
     \gamma_2=\e_2+\e_3  >\gamma_1=\e_2-\e_3. 
$$

The strongly compact and highest non-compact
roots
are constructed explicit in 
\cite[Table 4.7]{GW} for all quaternionic
Lie algebras, and
by checking the list we
have the following;
this can also be proved abstractly.

 \begin{lemm+}  \label{lem1-1}
   Let $\fg$ be an  irreducible
   quaternionic Lie algebra of rank $4$.
  The highest non-compact root  $\gamma$,
  the highest root $\beta_1$ and the simple non-compact root
$  \gamma_1$ are  related by   $\gamma=\beta_1-\gamma_1$. 
The roots
  \begin{equation}
  \label{gamma-s}
   {\beta_j}=
\gamma -\gamma_j =\beta_1 -\gamma_1-\gamma_j, \, 
2\le j\le 4,   
\end{equation}
form a system of strongly compact 
orthogonal roots
for the ideal $\fm_c^{\mathbb C}$ of $\fk^{\mathbb C}=\fsl(2, \mathbb C) + \fm_c^{\mathbb C}$.
The compact roots $\{\beta_j\}_{j=1}^4
$ are of the same length as
 $\{\gamma_j\}_{j=1}^4$
with  the inner products
$$
(\gamma_j, \beta_1)=(\gamma, \beta_1)=1,\, 1\le j\le 4,
$$
and
$$(\gamma, \gamma_1)=-1,    (\gamma, \gamma_j)=(\gamma, \beta_j)=1, \,
2\le j\le 4.$$
The following relations hold
  \begin{equation}
    \label{beta-gamma}
\beta_1=\frac 12(\gamma_1+\gamma_2+\gamma_3+\gamma_4),
\gamma=\frac 12(\beta_1+\beta_2+\beta_3+\beta_4).
  \end{equation}
In particular the  co-roots ${\gamma^\vee}$  is obtained as
  \begin{equation}
    \label{H-def}
    H:= {\gamma^\vee}=\frac 12 ({\beta_1^\vee}+
    {\beta_2^\vee}+
    {\beta_3^\vee}+{\beta_4^\vee}).    
  \end{equation}
\end{lemm+}

\subsection{Split real Lie algebras $\fg=\fe_{6(6)},
  \fe_{7(7)},
  \fe_{8(8)}  $}

   The real root system
  of $(\fg, \fk)$ is the same
  as the complexifed Lie algebra
  $\fg^{\mathbb C }=\fe_{6},
  \fe_{7},
  \fe_{8}  $ and  is  somewhat
  easier.
  So we fix $\gamma$ the
  highest non-compact root and the root
  vector $E_{\pm \gamma}$ as above.
There exists also
a system $\{
  \beta_1,\beta_2, \beta_3, \beta_4\}$
  of strongly orthogonal compact
  roots. Some precise construction of
  $\gamma$ and 
  $\{
  \beta_1,\beta_2, \beta_3, \beta_4\}$
is done  in  \cite{BK},
and there are Harish-Chandra
strongly orthogonal roots
for the Hermitian symmetric pair $(\fk,
\fl_1)$ defined 
in    (\ref{def-chi})  below,
and the relation
$\gamma=\frac 12(\beta_1+\beta_2+\beta_3+\beta_4)$
above still holds \cite[p. 2471]{BK}.

\section{Heisenberg parabolic algebras and groups}
\label{h-section}
\subsection{Root space decomposition of $E=E_\gamma + E_{-\gamma}$}
Recall that 
$\gamma$
is the {\it   highest
  non-compact  root} with the corresponding real form
$\fsl(2, \mathbb R)_\gamma=\mathbb R E +\mathbb R F +\mathbb R(i H)$ in   (\ref{X0-Y0}).     Let  $     \fa=\mathbb R E\subset \fp$ and consider
       the root space decomposition of $\fg$ with respect to $\fa$,
       \begin{equation}\label{heisen-par}
    \fg=\fn_{-}\oplus (\fm\oplus\fa) \oplus \fn
           \end{equation}
               where $\fn$ is the sum of  positive root spaces 
              and $\fm \oplus \fa$ is the zero eigenspace of $E$. 
              The Lie algebra  $\fn=\fn_2 +\fn_1$
              is a Heisenberg algebra \cite{GW}
              and $\fm +\fa+\fn$ is also called Heisenberg parabolic.
              We identify  the dual space $(\fa^{\mathbb C})^\ast$ with 
              $\mathbb C$ via $\nu \in \mathbb C \to \nu E^\ast: E\mapsto \nu$.

              The following formulas
              are obtained by taking $\fp$
              and $\fk$-parts
              of the decomposition     $$\fg= (\fn_2 +\fn_{-2}) +
              (\fn_1 +\fn_{-1}) + (\fm +\fa).$$

            \begin{lemm+} \label{p012-k012}
                There exist  orthogonal decompositions of $\fp$
  and $\fk$ as
                            \begin{equation}  \label{pk-dec}
  \fp=\fp_2\oplus\fp_1\oplus \fp_0,\, 
                \fk=\fk_2\oplus\fk_1\oplus \fl, 
              \end{equation}
              with 
              \begin{equation}
                  \label{k2}
              \fp_2=\mathbb RE +\mathbb RF,\,   \fk_2=\mathbb RiH, 
              \end{equation}
              \begin{equation}
    \fp_0=\fp\cap \fm,   \quad \fl=\fk\cap \fm,
\end{equation}  
such that the root spaces
$\fn_{\pm 1}
$ and $\fn_{\pm 2}
$
are given by
\begin{equation}
  \label{n12}
\fn_{\pm 1}=\{X \pm  [E, X]; X\in \fp_1\}, \quad
\fn_{\pm 2}=\mathbb R X_{\pm 2}, \quad X_{\pm 2}:=
\frac 12 (F\mp iH).
\end{equation}  
Moreover
$X\in \fp_1\mapsto  [E, X]\in \fk_1$  
is an isometry from Euclidean
space $(\fp_1, B(\cdot, \cdot))$
onto  $(\fk_1, -B(\cdot, \cdot))$.          \end{lemm+}              

The half sum    of positive roots is by definition
              \begin{equation}
                \label{rho-g}
                \rho_{\fg} = 1 +\frac 12 \dim \fn_1. 
              \end{equation}

\subsection{Heisenberg parabolic subgroup $P=MAN$}

Recall our  realization of 
group $G$ as $G=\Ad(\fg)$.
Let $A:=\{e^{tE}; t\in \mathbb R\}=
             \mathbb R^+$. The subgroup $\{g\in G; g E:=\Ad(g) E=E\}$
is a product $MA$, and we define also
           \begin{equation}
             \label{p210}
             P:=MAN, \,  L:=M\cap K=Z_{K}(\fa)
             =\{k\in K; \Ad(k) E\}. 
\end{equation}
Clearly their Lie algebras are $\fa$, $\fm$, $\fm +\fa+\fn$, and
$\fl=\fm\cap \fk$, respectively.
We denote 
\begin{equation}
  \label{u1-H}
U(1):=\exp(\mathbb RiH).
\end{equation}
The Lie group with Lie algebra $\mathbb R iH\subset \fk$
and this copy of $U(1)$ will be used
throughout the rest of our paper.

Let $L_0\subset L$, $M_0\subset M$ and
$M_0AN\subset P=MAN$
be the connected component of $L$, $M$ and 
$P$, respectively. In \cite[Lemma 4.2.2]{F} an explicit element $w_0\in 
      \exp(\fk)$
      (the element $w_1^2$ or $w_2^2$ there, representing 
the same class in $L/L_0$) 
 is constructed such that $\Ad(w_0)E=E, 
  \Ad(w_0)F=-F$ for any non-Hermitian 
   Lie algebra $\fg$.  We formulate
   this fact as follows.
   \begin{lemm+}
    The group $M$ has two components,  $M=M_0\rtimes \mathbb Z_2$. 
  So is also the group $L$, $L=L_0\rtimes \mathbb Z_2$ and
\begin{equation}\label{L-0-des}
  L_0=\{k\in K; k E=E, kF=F\}=
  \{k\in K; k E_{\pm\gamma}=E_{\pm \gamma}\}=M_0\cap K.
\end{equation}
\end{lemm+}
   \begin{proof}
     By definition both groups $M$  and $M_0$
     act 
     on $\fm$ by the adjoint representation
     with the subalgebra $\fm$ itself as a Lie algebra and
     $Ad(M)|_{\fm}= Ad(M_0)|_{\fm}$
     is
the adjoint group  of $\fm$. The group $M$ acts
     on $\fn_2=\mathbb RX_2$ by a character $\sgn g=\pm 1$,
     $\Ad(g) (X_{2})=\sgn(g) X_{2}$. In particular
     the connected component $M_0\subset \Ker \sgn\subset M$
     is a normal subgroup. Conversely if $g\in \Ker \sgn$ then
     $\Ad(g)|_{\fm} = \Ad(h)|_{\fm}$ for some $h\in M_0$. We claim that $g=h$,
     i.e. $h^{-1}g=1$ as element in $\Ad(\fg)$. 
 We have  $\Ad(h^{-1}g)|_{\fm}=1$,
     $\Ad(h^{-1}g)|_{\fn_2}=1$.
     The element    $\Ad(h^{-1}g)$ acts then trivially
     on $X_2$. But the adjoint action of 
     $Ad(M_0)|_{\fn_1}$    on $\fn_1$ is defined by
     the Lie algebra action of $\fm$ on $\fn_1$, thus
     $\Ad(h^{-1}g)|_{\fn_1} =1$ also;
     similarly
     $\Ad(h^{-1}g)|_{\fn_{-1}}
     =1, \Ad(h^{-1}g)|_{\fn_{-2}} =1
     $.  Now
     $\fg=(\fn_1 + \fn_2) + (\fm +\fa)
     +(\fn_{-1} + \fn_{-2})$
     and we find $\Ad(h^{-1}g)=1$ on $\fg$
     and $g=h\in M_0$.
     This proves $M/M_0=\mathbb Z_2$.
     It follows also from
\cite[Lemma 4.2.2]{F}
that $w_0$ does not commute with $M_0$ hence  $M=M_0\rtimes \mathbb Z_2$
as a semi-direct product of $M_0$ and $\mathbb Z_2=\{1, w_0\}$.
 
 The same argument proves that  $L=L_0\rtimes \mathbb Z_2$ and
    the connected component $L_0$ is given by
    $  L_0=\{k\in K; k X_2=X_2\}$. Since $K$
    acts on $\fp$ and $\fp^{\mathbb C}$, 
 the isotropic subgroup $L_0$
of  $E,  F\in
 \fp$
is the same as that of $X_{\pm 2}\in \fn_{\pm 2}$, 
and further the isotropic subgroup
of 
 $ E_{ \pm \gamma} 
\subset
\fp^{\mathbb C}$.
This completes the proof.
\end{proof}

\subsection{The subgroup $L_1\subset K$ and  character $\chi$}
We introduce now the subgroup $L_1$ and
the space $L^2(K/L_1; \chi^p) $.
\begin{defi+}
  \label{L-2-spaces}
  \begin{enumerate}
\item  We define the  subgroup     $L_1\subset K$
  and a character $\chi: L_1\to  \mathbb C$ by
  \begin{equation}
    \label{def-chi}  
    L_1:=\{h\in K; \exists \chi(h)\in \mathbb C, hE_\gamma=\chi(h) 
    E_\gamma \}.
  \end{equation}
\item Denote $L^2(K/L_1; \chi^p) $
  the $L^2$-space of sections of
  the line bundle $K\times_{(L_1, \chi^{p})}\mathbb C$ 
over $K/L_1$,  namely it is space of 
 $L^2$-functions $f$ on $K$ 
such that 
$$
\chi^{p}(h)f(kh) = f(k), \quad k\in K, h\in L_1. 
$$
\end{enumerate}  
   \end{defi+}  
   It is obvious that $\chi: L_1\to U(1)\subset \mathbb C$ is
   a unitary character.
   It follows by definition $L_0=\ker \chi\subset L_1$ is a normal
   subgroup of $L_1$, the subgroup  $L_1$ is  the stabilizer  in $K$
   of the line $\mathbb CE_{\mathbb C}$
   in the projective space $\mathbb P(\fp^{\mathbb C})$,
   $L_0$ is the stabilizer of the point $E_\gamma$,
   and $L$    the point $E\in \fp$.
   Namely
   $$K/L\subset S(\fp),
   \, K/L_0\subset \mathbb P(\fp), \,
   K/L_1\subset \mathbb P(\fp^\mathbb C),
   $$
as submanifolds of the  unit sphere $S(\fp)$,
  the real projective space and
  the complex projective space, respectively.

The Lie  algebra  $\fl_1$ is 
       $$
   \fl_1:=\{X\in \fk; \ad(X)E_\gamma=
    d\chi(X)E_\gamma \}
    $$
    and $\ad(X)E_\gamma=   d\chi(X)E_\gamma$ for $X\in \fl_1$.
    Here         $d\chi$ is the differential of the character 
    $\chi$.

    \begin{lemm+}
      \label{lem2.4}
      The element  $w_0\in K$ normalizes 
      $L_1$,    $\Ad(w_0)L_1=L_1$, 
     and the character $\chi$
     on $L_1$ satisfies $\chi(\Ad(w_0)k)=
     \chi(k)^{-1}$, $k\in L_1$.
     Moreover the following relations hold 
     $$
     L_0= L\cap L_1,\,      \, L_1={U(1)}\times L_0 
     =U(1) L_0, \, L_1/L_0={U(1)}, \fl_1=\fk_2+\fk_0
     $$
     where $\fk=\fk_2+\fk_1+\fk_0$  is as
     in Lemma \ref{p012-k012}.
\end{lemm+}      
\begin{proof}
  The element $w_0$ satisfies $w_0 E=E, 
  w_0 F=-F$ and then 
  $w_0 E_{\pm \gamma}=E_{\mp\gamma}$, so the 
 normalizer properties  of $w_0$ follow by definitions. 
 We prove then $L_0=L\cap L_1$. It is clear that $L_0\subset L\cap
 L_1$.
 Suppose  $k\in L_1\cap L$. 
   Then 
   $kE= E$ and $k$ preserves the $\pm 2$-root spaces 
   $\fn_{\pm 2}=\mathbb R(F\mp iH)$ of $\ad(E)$. 
   But $k$ preserves also $\fp$ and $\fk$
   and is an orthogonal transformation, so $k F=\pm F$. 
   If $k F=-F$ then $kE_{\gamma}=
   \frac 12 
  k(E -i F)= \frac 12(  k(E) -i k(F))=
   \frac 12 ( E +i F)= E_{-\gamma}$, which is not 
   a complex scalar multiple of $E_{\gamma}$ and is contradiction 
   to  $k\in L_1$. Thus $k F=F$ and $k E_{\gamma}
   = E_{\gamma}$ and $k\in L_0$.   This proves the equality $L_0=L\cap L_1$.
   Finally the character $\chi$   defines a natural isomorphism $L_1/L_0=U(1)$,
   which implies also $\fl_1=\fk_2+\fl$ since
   $\fk_2=\mathbb i RH$.
   \end{proof}

\subsection{Compact Hermitian 
  Symmetric Space $K/L_1$}

We recall the following known facts;
see \cite{Ch, F, Sl-St}
for the general cases
of Heisenberg Lie algebras
in simple Lie algebras, 
\cite{GW} for the quaternionic
Lie algebras and 
 \cite{BK} for exceptional 
 Lie algebras $\fe_{6(6)},
 \fe_{7(7)}, 
 \fe_{8(8)}
 $.
Recall that  $iH= i{\gamma^\vee}$. 
Let
${U(1)_{\beta_1}}=
\exp(\mathbb Ri H_{\beta_1})\times 
1/((-1)_{\beta_1}\times \epsilon)\subset K.$

\begin{prop+}
  \label{prop-1} Let $G$
  and $\fg$ be as in the Table
  \ref{tab:1}. 
We have  the following commutative
diagram of  $K$-equivariant coverings 
    of homogeneous spaces of $K$, 
      $$
\begin{tikzcd}[arrows=->]
  {G/M_0AN= K/L_0} \arrow{r}{\mathbb Z_2}
  \arrow{d}{U(1)}
  & {G/P=G/MAN=K/L}
  \arrow{d}{U(1)}
\\
K/L_1 
=K/U(1)L_0 
\arrow{r}{\mathbb Z_2}
&
K/{U(1)}L 
\end{tikzcd}. 
$$
The space   $K/L_1=K/{U(1)}L_0$
is a compact Hermitian symmetric space 
of rank $4$
with Cartan involution
$\Ad(\exp(\pi iH))$
and decomposition  $\fk=(\fk_2 + \fl) +\fk_1$,
the complex  structure 
$\ad(iH)$
on $T_e(K/L_1) =\fk_1$, 
and Harish-Chandra
 strongly orthogonal roots
 $\beta_1, \beta_2, \beta_3, \beta_4$ constructed above. 
 When $G/K$ is quaternionic
the space  $K/L_1=K/{U(1)}L_0$
biholomorphic to  the complex manifold 
$\mathbb P^1\times M_c/L_c$, 
 and $M_c/L_c$ is a compact 
symmetric space of $M_c$ of rank 3 (whose 
non-compact dual is of  tube type).
\end{prop+}
\begin{proof} 
 It is a general fact that $G/P=G/MAN$ as a $K$-space
  is $G/P=K/L$. Now $M_0\subset M=M_0\rtimes \mathbb Z_2$
  and $L_0\subset L=L_0\rtimes \mathbb Z_2$
  by Lemma \ref{lem2.4},
  thus $G/M_0AN=K/L_0\to  G/M=K/L$ is a $\mathbb Z_2$-covering.
  It follows also from definitions that $K/L_0\to K/L_1$ is a circle bundle. Clearly
  these covering are $K$-equivariant and we get the commutative diagram.
It is  known that $\fg $ has the subalgebra
 $\fsl(2, \mathbb R)_{\gamma} +\fm$ as a non-compact symmetric
subalgebra; see e.g. \cite[Corollary 2,9]{Sl-St} and \cite{F}.
Restricted to $\fk$ we get a symmetric
subalgebra $\mathbb R i H + \fl
=\fk_20\fl$ 
 of $\fk=(\fk_2 +\fk_0)+\fk_1$. 
Thus $K/U(1)L_0$ is a symmetric space.
The compact strongly orthogonal compact
roots are $\beta_1,   \beta_2, 
 \beta_3, \beta_4$; when $\fg$
 is quaternionic the roots
  $ \beta_2, 
 \beta_3, \beta_4$ are compact
strongly orthogonal roots for $\fm_c$
by the definition of $\fk$,
and $\fm_c$ has a compact symmetric pair
$(\fm_c, \fl_c)$ of rank 3
of tube type
where $\fl_c=\fm_c\cap (\mathbb RiH +\fl)$; more precisely
the non-compact dual of
$(\fm_c, \fl_c)$ is of tube type;  see \cite[Lemma 6.9]{GW}.
It is also a  general fact (\cite[Chapter VIII]{He},  \cite{Loos},
\cite[p.221]{Up})
that 
 the Cartan involution and
 complex structure of an
irreducible Hermitian symmetric compact or non-compact
  pair of tube-type 
  are determined by the Harish-Chandra strongly orthogonal
  roots, in our case $\beta_1, \cdots, \beta_4$,
 by
    $\Ad(\exp(\frac{\pi}2(
    {\beta_1^\vee} +
    \cdots+    {\beta_r^\vee}))
    $
and by
    $\frac i2 \ad ({\beta_1^\vee} +
    \cdots+    {\beta_r^\vee})$. But
 $\frac i2 (\beta_1^\vee +
  \beta_2^\vee +
  \beta_3^\vee +
  \beta_4^\vee ) = \frac 12 i H$, by  (\ref{H-def}).
  When $\fg$ is quaternionic
  the Lie algebras $\fsu(2)_{\beta_1}$ and $\fm_c$ are commuting  so we have
  $K/L_1=K/{U(1)}L_0$
biholomorphic to  the complex manifold 
$\mathbb P^1\times M_c/L_c$
with $\mathbb P^1=SU(2)_{\beta_1}/U(1)_{\beta_1}$.
\end{proof}  
\begin{rema+}
  We note that for quaternionic $\fg$, 
  $M/L=M_0/L_0$ is a
  non-compact symmetric pair dual to
  $M_c/L_c$. However the subgroups
  $L$  and $L_c$ do not coincide in our
  realization as subgroups of $G$.
 Using the action $\ad iH$ 
we get the Cartan decomposition of the
symmetric pair $(\fk, \fl_1)$,
$$
\fk =\fsu(2) + \fm_c
=\fl_1 +\fq, \,
\fsu(2) = \mathbb Ri\beta^\vee \oplus + \fq_{\fsu(2)},
\fm_c=\fl_c +\fq_c,
\fq=\fq_{\fsu(2)} +\fq_c, 
$$
and 
a refined decomposition
of $\fn_1$, 
$\fn_1=(\fn_{2, -1} + \fn_{-1, 2}  ) + \fn_1'$
where
$$
\fq_{\fsu(2)}=\{X +\theta (X); X\in \fn_{2, -1} +
\fn_{1, -2} 
\}, 
\fq_c=\{X +\theta (X); X\in \oplus \fn_{1}'\}.
$$
See  \cite[Section 13]{GW} and \cite[Section 4.2]{F}
for the details.
\end{rema+}

\subsection{$\rho_\fg$ and $\rho=\rho_{\fk^{\mathbb C}}$}
     Let $\rho_{\fk^{\mathbb C}}$
     be the half-sum  
     of the positive roots of $\ft^{\mathbb C}$ in $\fk^{\mathbb C}$.
     If $\fg$ is quaternionic and $\fg\ne \fso(4, d)$ then
    $ \rho_{\fk^{\mathbb C}}$
    can be computed directly from the 
Table \ref{tab:1}    (see also \cite[p.229]{Up}),
  \begin{equation}
    \label{rho-k}
    \rho_{\fk^{\mathbb C}}
    =\frac 12\sum_{j=1}^4\rho_j\beta_j,\,
\rho_1=1, \rho_2=1+2a,  \rho_3=1+a, \rho_4=1;
\end{equation}
where $a$ is
the multiplicity for the (restricted) root
$\frac{\beta_i-\beta_j}2, 2\le i\ne j\le 4$ and
a list of the root multiplicity $a$ is found in
the Table \ref{tab:2}.

Let $\fg= \fso(4, d)$. We follow the notation
in Section \ref{strong-ort}.
The half sum of positive roots
  $$
  \rho=\rho_{\fk^{\mathbb C}} =
  \frac 12(\beta_1 +\beta_2 
  +  (d-3)\beta_3 +\beta_4), 
  $$
with
\begin{equation} 
  \label{rho-k-so}
  \rho_1=1, \, \rho_2= 1, 
\rho_3= d-3, \, 
\rho_4=1.
\end{equation} 

For $\fg=\fe_{6(6)},
\fe_{7(7)},
\fe_{8(8)},$
we have 
  \begin{equation}
    \label{rho-k}
    \rho_{\fk^{\mathbb C}}
    =\frac 12\sum_{j=1}^4\rho_j\beta_j,\,
\rho_1=1+3a, \rho_2=1+2a, \rho_3=1+a, \rho_4=1.
\end{equation}

We shall need the following relation
between the half sums
$\rho_\fg$ and $\rho=\rho_{\fk^{\mathbb C}}$, which
can also  be checked easily from the Table \ref{tab:2}.
\begin{lemm+}     \label{rho-rho}
 We have 
  \begin{equation}
  d=  \frac 12 \dim\fn_1 
  = \sum_{j=1}^4\rho_j, \,
  \rho_{\fg}= 1+\frac 12 \dim\fn_1=  1+\sum_{j=1}^4\rho_j.
  \end{equation}
\end{lemm+}  
  \begin{proof} We compute  $\tr\ad H$ on $\fk_+^{\mathbb C}$. 
  It follows by definition that $\tr\ad H= 2\rho_{\fk^{\mathbb C}} (H)$. 
  But by Proposition \ref{prop-1}
  $\ad iH$ defines also the complex structure 
  of the Hermitian symmetric pair $(\fk, \fl_1)$
  so that $
  2\rho_{\fk^{\mathbb C}} (H) 
  =\tr\ad H_{\fk^{\mathbb C}_+}= \dim_{\mathbb C} K/L_1 
  =\frac 12(  \dim_{\mathbb R} K/L -1) =
  \frac 12(\dim_{\mathbb R} \fn -1) 
  =\frac 12(\dim {\fn_1}) =\rho_{\fg}-1$. 
  Thus  $ \rho_{\fg} =   2\rho_{\fk^{\mathbb C}} (H)+1$. 
  When $\rank \fg=4$ all roots $\beta_j(H)=1, 1\le j\le 4$
  and we find $\rho_{\fg}= 1+ 2\rho_{\fk^{\mathbb C}}(H) = 1+
  \rho_1 +\rho_2+\rho_3+\rho_4$ by (\ref{rho-k})
  and (\ref{rho-k-so}).
  This completes the proof.
  \end{proof}
  \section{Induced Representations
    $\pi_\nu= \Ind_P^G(1\otimes \exp(\nu)\otimes 1)$
    and its $K$-types
  }
\label{deco-sect}

The main purpose of 
this paper is to study 
the unitarity 
of induced  representations  $\Ind_P^G(1\times e^\nu \otimes 1)$
from the Heisenberg parabolic subgroup $P=MAN$
and their subrepresentations.
This will be done by  finding the 
parametrization for the $K$-types in 
the representations and  performing explicit 
computations on $\fa$-actions  
on $L$-invariant elements  of fixed  $K$-types.
Since the group $M$
 is not connected there is also another family
 $\Ind_P^G(\sgn \times e^\nu \otimes 1)$
 with the $\sgn$ representations of $M$ which
 we treat in the same time.

\subsection{Induced Representations 
of $G$ and their $K$-types
}
Recall \cite{Kn} that
the induced representation 
       $I(\nu)=(\pi_{\nu}, G)=\Ind_{P}^G(1\otimes e^\nu\otimes 1)$  of $G$
       is defined as
       left regular representation
 $\pi_\nu(g) f(h) = f(g^{-1}h), g, h\in G$, 
      on the space of measurable functions $f$ on $G$
such that 
$$
f(gm\exp(tE)n) = e^{-t\nu} f(g), \quad t\in \mathbb R, m\in M, \quad n\in N 
$$
and $f|_k\in L^2(K)$.  Likewise denote
$I(\nu, -)=(\pi_{\nu, -}, G):=\Ind_P^G(\sgn\otimes e^\nu\otimes 1)$
with $\sgn$ the  representation of $M$ 
via $M\to M/M_0=\mathbb Z_2=\{\pm 1\}$.  We abuse
notation and denote 
 the corresponding subspaces of $K$-finite elements, i.e.
the underlying Harish-Chandra modules also by $I(\nu) $ and  $I(\nu, -)$.

For computational convenience we introduce also
$I(\nu, 0)=(\pi_{\nu, 0}, G)=\Ind_{M_0AN}^G(1\otimes e^\nu\otimes 1)$
induced from the connected component $M_0$ of $M$.
The right regular action  of $w_0$, 
${w_0}: f(g)\mapsto {w_0}f(g)=f(gw_0)$
preserves the spaces
$\Ind_{P}^G(1\otimes e^\nu\otimes 1)$
and $\Ind_{P}^G(\sgn \otimes e^\nu\otimes 1)$
since $\Ad(w_0)(E)=E$,
$\Ad(w_0)(M)=M$, $\Ad(w_0)(N)=N$. We have thus
\begin{equation}
  \label{I-pm}
I(\nu, 0)
=I(\nu) +I(\nu,  -),
I(\nu, 0) =L^2(K/L_0), 
I(\nu) =L^2(K/L), \,
I(\nu,  -)=L^2(K/L; \sgn)
\end{equation}
as representations of $\fg^{\mathbb C}$
and $K$; the elements $f\in L^2(K/L)\subset L^2(K/L_0)$ will
be called even and $ f\in L^2(K/L_1)\subset L^2(K/L_0)$ odd.
Here
$L^2(K/L; \sgn)
$ is the $L^2$-space
of functions $f\in L^2(K)$ such
that $f(kh)=\sgn(h) f(k)$, $h\in L\subset M$.
Observe also that the right multiplication $R_{w_0}: f(k)\mapsto f(kw_0)$
defines a $K$-invariant isometric isomorphism
\begin{equation}
  \label{right-w}
R_{w_0}:  L^2(K/L_1; \chi^{p}) 
\to 
L^2(K/L_1; \chi^{-p}).
\end{equation}

We shall use the generalization
of the Cartan-Helgason theorem
by Schlichtkrull
to find the   decomposition for  $L^2$-spaces 
$L^2(K/L_1; \chi^p)$
of line bundles 
over $K/L_1$.
The result of Schlichtkrull
\cite{Sch}
is stated
using the Cartan subalgebra
for the non-compact dual of $K/L_1$.
To ease notation we use
an equivalent formula by taking the fixed
Cartan $\ft=\ft_{\fsu(2)} +\ft_c\subset \fk\subset \fg$,
with $\{\beta_j\}_{j=1}^4 $ the strongly
orthogonal compact roots.
We apply the result
in \cite{Sch} (see also \cite{Shi})
to our case and get

  \begin{lemm+}
    \label{C-H-S}
    \begin{enumerate}
      \item
As $K$-representation spaces we have
    \begin{equation}
      \label{L2KL+-0}
      L^2(K/L_0)=
      \sum_{p=-\infty}^{\infty}
      L^2(K/L_1, \chi_p),
          L^2(K/L_1, \chi_p)=
      \sum_{p, \mu} W_{\mu; p}, 
  \end{equation}
      \begin{equation}
      L^2(K/L)=\sum_{p\ge 0, \mu} W_{\mu; p}^{+}, 
    L^2(K/L; \sgn)=
    \sum_{p>0, \mu} W_{\mu; p}^{-},
  \end{equation}
  where
  $W_{\mu, p}^{\pm}$
  is obtained
  as the even and odd part in 
  $W_{\mu, p} +W_{\mu, -p}  $
  by our convention above,  and 
  $W_{\mu, p}$ is of highest weight
$$
\mu= (\mu_1, \cdots, \mu_4):=
\frac 12 (\mu_1 
\beta_1 
+\mu_2 \beta_2 
+\mu_3 \beta_3 
+\mu_4 \beta_4), 
$$
subject to the following conditions:
\begin{enumerate}
  \item
If $\fg$ is quaternionic and $\fg\neq \fso(4, d)$,
\begin{equation}
  \label{rank4}
    \mu_1 \ge |p|, \,  \mu_2\ge \mu_3 \ge 
\mu_4\ge  |p|, \,  \mu_j=p \pmod 2 \forall j.  
\end{equation}
  \item $\fg=\fso(4, d)$,  $d>4$, 
\begin{equation}
  \label{so-4-d} 
\mu_1,  \mu_2\ge |p|, \, 
 \mu_3\ge   \mu_4\ge |p|, \, 
\mu_j=p \pmod 2 \forall j.
\end{equation}
If  $d=4$,  $\fk^{\mathbb C}=
\fsl(2, \mathbb C) _{\beta_1}
+
\fsl(2, \mathbb C)_{\beta_2} +
\fsl(2, \mathbb C)_{\beta_3} +\fsl(2, \mathbb C)_{\beta_4}$, 
\begin{equation}
\label{so-4-4}   
\mu_j \ge |p|, 
\mu_j=p  \pmod 2  \forall j.
 \end{equation}
\item If $\fg=\fe_{6(6)},
  \fe_{7(7)},  \fe_{8(8)}$,
\begin{equation}
  \label{rank4-exept}
    \mu_1 \ge \mu_2\ge \mu_3 \ge 
\mu_4\ge  |p|, \,  \mu_j=p \pmod 2 \forall j.  
\end{equation}  
\end{enumerate}

\item For each integer $l, l=p \pmod 2$
  satisfying the same condition as $p$ above 
  the space  $W_{\mu; p}$ contains,
  up to non-zero scalars,  a 
  unique nonzero element  $\phi_{\mu; l, p}$
    such that
  $$
  \chi^{l}(h) \chi^{p}(h') 
  \phi_{\mu; l, p}(hkh') 
  =
  \phi_{\mu; l, p}(k), h, h'\in L_1, k\in K.
  $$
  The even and odd parts
  of $\phi_{\mu; l, p}$ are
  \begin{equation}
          \phi_{\mu; l, p}^+
:=\phi_{\mu; l, p} + {w_0} \phi_{\mu; l, p}\in 
     W_{\mu; p}^+, p\ge 0; \, 
     \phi^-_{\mu; l, p}:=\phi_{\mu; l, p} - {w_0} \phi_{\mu; l, p}
          \in W_{\mu; p}^-, p>0. 
   \end{equation}
\end{enumerate}
 \end{lemm+}

We shall need in particular the spherical polynomials
in the space  $\fp^{\mathbb C}$:
It is of highest weight 
$\gamma=\frac 12(\beta_1+\beta_2+\beta_3+\beta_4):=1^4$
with highest weight vector $E_\gamma$, 
as representation of $K$ with $(E_\gamma, E_{-\gamma})=1$. The 
$L_1$ action on $E_{\pm\gamma}$ satisfies 
$$
\chi(k)  k^{-1} E_{\gamma}
=E_\gamma, \, \chi(k)^{-
  1}  k^{-1} E_{-\gamma}
=E_{-\gamma}, \, k\in L_1,$$
by our definition. 
Hence
the matrix coefficient
$\psi_{1^4}(k):= 
\langle k^{-1}E, E\rangle$ is
\begin{equation}
  \label{eq:1-4-sph}
\psi_{1^4}=\( \phi_{1^4; 1, 1} +
  \phi_{1^4; 1, -1}\) +
  \( \phi_{1^4; -1, 1} +
   \phi_{1^4; -1, -1}\) 
   ={\phi_{1^4; 1, 1}^+}  +   {\phi_{1^4; -1, 1}^+},
  \end{equation}
  is a sum of $4$ spherical polynomials.
 The matrix coefficient
  $  \langle k^{-1}E, F\rangle $ is 
\begin{equation}
  \label{eq:1-4-sph-}
  \begin{split}
    \langle k^{-1}E, F\rangle
    &     = \langle k^{-1}(E_\gamma + E_{-\gamma}),
   i( E_\gamma -E_{-\gamma})    \rangle
    \\
    &=
-i  \(
  \phi_{1^4; 1, 1} (k)-
  \phi_{1^4; 1, -1}(k)
  \) -i
  \( \phi_{1^4; -1, 1}(k) -  \phi_{1^4; -1, -1}(k)\) \\
  &
 : =-i   \phi^-_{1^4; 1, 1}(k)  -i   \phi^-_{1^4; -1, 1}(k),
\end{split}
\end{equation}
 a sum of 
 ${\phi^-_{1^4; \pm 1, 1}}\in W_{1^4; \pm 1, 1}^-$.

\subsection{Casimir elements}
\label{cas-sec}

We shall use the technique  
developed in \cite{B-O-O}
to treat certain differentiation
of spherical polynomials.
The technique can roughly be summarized
as first degree
differentiations
can be computed using certain
second order Casimir elements.
To introduce these elements we  
recall first the  decompositions (\ref{p210})  
and the correspondence between $\fp_j$,  $\fk_j$ 
and  $\fn_j$,  $j=0, 1, 2$, described in Lemma \ref{p012-k012}.  
So let $\{u_j\}\subset \fp_1$ 
be an orthonormal basis (ONB)  
of $\fp_1$ with respect to the Killing form  
and let $u_j+z_j\in \fn_1$ 
be the corresponding basis of $\fn_1$, $(u_i, u_j)=\delta_{ij}$.  
Then $z_j=\ad (E)u_j$, $u_j=\ad (E) z_j=\ad^2(E) u_j$  and $\{z_j\}$
is an ONB
of $\fk_1\subset \fk$ 
in the sense that $(z_i, z_j)=-\delta_{ij}$.

\begin{defi+}
  \label{defi-cas}
  Let $\fk=\fk_2+\fk_1+\fl$ be  
as in Lemma \ref{p012-k012}. 
 Define the Casimir element of the subspace $\fk_1$ by 
$$
\Cas_1=\sum_j z_j^2, 
$$ 
where $\{z_j\}$ is any ONB for $\fk_1$. 
\end{defi+}

 Denote  $\Cas_2$ and $\Cas_0$  the Casimir element 
 of the Lie algebras $\fk_2=\mathbb Ri H$
 and $\fl$   respectively. 
 Then $\Cas_2=-\frac 12 H^2$ since $\fk_2=\mathbb RiH$,  
and $(H, H)=2$, and the Casimir element of $\fk$ is 
$$
\Cas:=\Cas_{\fk}=\Cas_2+\Cas_1+\Cas_0. 
$$
The eigenvalue 
$\Cas(\mu)$ of $\Cas$ on $W_\mu$ is
\begin{equation}
  \label{cas-mu}
\Cas(\mu): = -(\mu + 2\rho, \mu) 
=-\frac{1}2\sum_{j=1}^4 (\mu_j +2\rho_j)\mu_j,
\end{equation}
since $\mu=\frac 12 (\mu_1 \beta_1 +
\cdots + \mu_4 \beta_4)$, with the Killing form 
$(\beta_j, \beta_j)=2$ and $\{\beta_j\}_{j=1}^4$
strongly orthogonal.

 \section{Lie algebra action of $\fg$
   on the induced representations}
 \label{main-sect}

 \subsection{ The action of    $E  $
   on  $I(\nu)$ and   $I(\nu, -)$
 }
We compute now the action of $E$ on  $L$-invariant elements
$ {\phi^{\pm}_{\mu; l, p}}$, which shifts
the parameter 
${(\mu; l, p)}$ to
${(\mu+\sigma; l+\delta_1, p+\delta_2)}$, $\sigma=(\sigma_j)_{j=1}^4\in \mathbb Z_2^4=\{\pm 1\}^4, \delta_1, \delta_2=\pm 1$.
We denote 
\begin{equation}
  \label{cri-coe}
  A_{\nu; \mu; l,  p}(\mu+\sigma;  l+\delta_1, p+\delta_2)
  =\frac 12\left( \nu +\sum_{j=1}^4 
  (\mu_j +\rho_j)\sigma_j -
  \sum_{j=1}^4\rho_j  -2p\delta_2 
  \right)
  C_{\mu, p}(\mu+\sigma;  l+\delta_1, p+\delta_2)  
\end{equation}
where   $C_{\mu, p}(\mu+\sigma;  l+\delta_1, p+\delta_2)$
are the coefficients in Section \ref{appen-sph-lp},
Proposition \ref{prod-exp-phi}, 
which depend only on $|l|, |l+\delta_1|, |p|, |p+\delta_2|$.

 \begin{theo+}
   \label{pi-E-W}
       The action of $\pi_{\nu}(E) $ on 
    $f^{+}={\phi^+_{\mu; l, p}}\in W_{\mu, p}^{+}\in I(\nu) $, $p\ge 0$,
  is given by 
  \begin{equation}
    \label{cri-coe-1}
\pi_\nu(E) 
f 
=\sum_{\sigma\in \mathbb Z_2^4, \delta=(\delta_j)\in \mathbb Z_2^2}
A_{\nu; \mu, l, p}
(\mu+\sigma;  l+\delta_1, p+\delta_2) 
{\phi^+_{\mu+\sigma; l+\delta_1, p+\delta_2}}.
  \end{equation}
   The action of $\pi_{\nu, -}(E) $ on 
    $f^-={\phi^-_{\mu; l, p}}\in W_{\mu, p}^{-}\in I(\nu, -) $, $p>0$,
  is given by 
\begin{equation}
  \label{cri-coe-2}
  \pi_{\nu, -}(E) 
f^- 
=\sum_{\sigma\in \mathbb Z_2^4, \delta\in \mathbb Z_2^2}
A_{\nu; \mu; l, p}
(\mu+\sigma;  l+\delta_1, p+\delta_2) 
{\phi_{\mu+\sigma; l+\delta_1, p+\delta_2}^{-}}. 
\end{equation}
\end{theo+} 

Here it is understood that
the terms $\phi^{\pm}_{\mu+\sigma; l+\delta_1, p+\delta_2}$
will not appear if
$({\mu+\sigma; l+\delta_1, p+\delta_2})$ does not satisfy
the conditions
(\ref{rank4})-(\ref{so-4-d})-(\ref{so-4-4})
and  $
\phi^{-}_{\mu; l, p} =0$ for $p=0$.

First of all it  follows by 
Lemma \ref{C-H-S} and the proof of Proposition  \ref{prod-exp-phi}
in Appendix
that  $\pi_{\nu}(E) f$
is a linear combination of
${\phi^{\pm}_{\mu+\sigma; l+\delta_1, p+\delta_2}}$
for $f=\phi_{\mu; l, p}^{\pm}$. This fact
will be used implicitly below.

\begin{rema+} By abstract argument
  we can prove that coefficient is
  an affine function of $\nu$ and
up to an non-zero constant it is
  of the form $\nu + c$. Observe also
  that $  \pi_{\nu, -}(E) f$
  involves {\it left differentiation} of $f$ as a function
  on $G$ restricted to $f$. The restriction to $K$
  is performed after the differentiation.
\end{rema+}

\begin{rema+}
  It follows from our theorem 
that the coefficients 
$A_{\nu; \mu; l, p}(\mu+\sigma;  l+\delta_1, p+\delta_2) $
 satisfy the following 
symmetry relation 
$$
A_{\nu; \mu+\sigma; l+\delta_1, p+\delta_2}
(\mu;  l, p) = -
A_{2\rho_{\fg}- \nu; \mu; l, p}(\mu+\sigma;  l+\delta_1, p+\delta_2). 
$$
This is also a consequence of a general 
fact that $\pi_{\nu}$ is unitary 
for $\nu=\rho_{\fg} + \tilde\nu, \tilde\nu\in i\mathbb R$. 
\end{rema+}

The rest of this section is 
devoted to the proof of the theorem. 
To simplify notation 
we write $\pi_{\nu, 0}(E), 
\pi_{\nu, \pm}(E) $ all as 
 $\pi(E)$ and we  
drop the sub-index $(\mu; l, p)$ and write  temporarily 
$$
C_{\sigma; \delta}: 
=C_{\mu + \sigma; l+\delta_1, p + \delta_2}, 
$$
for $\sigma\in \mathbb 
Z_2^4, \delta\in \mathbb Z_2^2$, 
where $C_{\mu + \sigma; l+\delta_1, p + \delta_2}$
is the coefficient in 
Proposition \ref{prod-exp-phi};
as 
$(\mu; l, p)$ is fixed so there will be 
no confusion. Note that 
$C_{\sigma; \delta}$ depends 
only on $|l|, |p|, |l+\delta_1|, |p + \delta_2|$
even though it is not reflected by the 
abbreviation. 
The  proposition 
 now reads as 
\begin{equation}
\label{prod-exp}
    \phi_{1^4; \delta_1, \delta_2}(k) 
    \phi_{\mu; l, p}(k) 
    =  \sum_{\sigma\in Z_2^4} C_{\sig; \delta}
{ \phi_{\mu +\sig; l +\delta_1, p+\delta_2}}(k)
\end{equation}
and the coefficients $C_{\sigma;  |l+\delta_1|, |p+\delta_2|} $
is non-zero whenever $(\mu + \sigma; l+\delta_1, p +  \delta_2)$
satisfies the  condition (\ref{rank4})-(\ref{so-4-d})-(\ref{so-4-4}).

We denote also the {\it right} differentiation of any element $X\in \fg$ on 
$f\in C^\infty(G)$  by $Xf$, namely 
\begin{equation}
\label{r-diff}
(X f)(k)=\frac{d}{dt}
f(g\exp(tX))|_{t=0}, g\in G. 
\end{equation}
 The differentiation can be extended
 to $X\in \fg^{\mathbb C}$ so 
 that $X\to Xf$ is complex linear. All differentiations
 will be right differentiation unless otherwise
 specified.

It follows by the definition of the induced representations that $Xf=0$ for any
$X\in\fm +\fn$ and $K$-finite
element $f\in I(\nu),  I(\nu, -)$
or $I(\nu, 0)$.
Each element $X\in \fp\subset \fg$, 
 is, mod $\fm+\fn$,  a sum of three elements
 in the  spaces $\fa, \fn_{-2}$ and  $\fn_{-1}$,
 respectively. Correspondingly  the action
 $\pi(E) {\phi_{\mu; l, p}}$
 is a sum of three terms $I, II, III$ (see below), the first two
 terms $I$ and $II$ require an expansion of the product
 of spherical polynomials, and the third
 term $III$ an expansion of the differentiations
 of spherical polynomials.

Any $K$-finite element $f\in L^2(K/L_1; \chi^p)$
transforms as $\chi(h)^pf(kh) =f(k)$,
thus
\begin{equation}
  \label{eq:H-on-f-0}
(iH)f(k) = -2ip f(k), \, Hf(k) = -2p f(k)  
\end{equation}
since $d\chi^p(H)= 2p $
as $[H, E_\gamma]=2 E_\gamma$.
Observe also that
\begin{equation}
    \label{eq:H-on-psi}
  (iH) \psi_{1^4}(k)=
\langle k^{-1} E, [iH, E]\rangle=
\langle k^{-1} E, 2i  
E_{\gamma}
- 
2i E_{-\gamma}
\rangle  
=2\langle k^{-1} E, F\rangle.
\end{equation}
Thus the eigenvalue $
\Cas_2=-\frac 12 H^2$ on the section $L^2(K/L_1; \chi^{\pm p})$
is 
\begin{equation}
  \label{cas-mu-2}
\Cas_2(p):= -2p^2.
\end{equation}
We shall also need
the Leibniz's rule, 
$$
\Cas_1 ( f g) 
= (\Cas_1 f ) g +
2\sum_j (z_j f)( z_j g) 
+f(\Cas_1  g ). 
$$
as differential operators.

\subsection{A general formula 
  for  $\pi(E)=\pi_{\nu, 0}(E)$}
 
Recall the decomposition  $\fp=\fp_2\oplus\fp_1\oplus\fp_0$
in Lemma  \ref{p012-k012}.
Let $P_2, P_1, P_0$ be the orthogonal projection
 onto the respective subspaces; the space
 $\fp_2=\mathbb R E \oplus
 \mathbb R F $, with
 $\langle
 E, E\rangle
 =\langle F, F
 \rangle=2$,
 and thus $P_2(X) =\frac 12 \langle X, E\rangle E +
 \frac 12 \lgl X, F\rgl F$, $X\in \fp$.

 \begin{lemm+} Each  element $X\in \fp$ has the following decomposition
   $$
   X =\frac 12 
   \lgl X, E\rgl  E + 
   \frac 12    \lgl X, F\rgl  (iH) +  [ P_1( X), E] \pmod{\fm+\fn}. 
   $$
\end{lemm+}
\begin{proof}    We have $P_0(X) \in \fm$,
  $F -i H\in \fn_2,   P_1( X) + [E, P_1(X)]\in \fn_1$,
  by   Lemma  \ref{p012-k012},
  i.e. $P_0(X)=0, F=iH,
   \, P_1(X)= -[E, P_1(X)] =[P_1(X), E]
    \pmod{\fm+\fn}$.
Therefore
\begin{equation*}
  X= P_2(X) + P_1(X)  =\frac 12 \lgl
  X, E\rgl   E + 
 \frac 12
\lgl X, F\rgl  F + 
   P_1( X) 
=\frac 12 
\lgl X, E\rgl   E + 
 \frac 12  \lgl X, F\rgl  (i H )
+ [P_1(X), E].
\end{equation*}
\end{proof}
For any $f\in I(\nu, 0) \subset C^\infty(G)$
we have then
\begin{equation}
\begin{split}  
  \label{Xf}
 \pi(X) f &=\frac 12 
   \lgl X, E\rgl  E f + 
   \frac 12    \lgl X, F\rgl  (iH)f +  [ P_1( X), E]f
   \\
     &=-\frac \nu 2   \lgl X, E\rgl  f + 
   \frac 12    \lgl X, F\rgl  (iH)f +  [ P_1( X), E]f
 \end{split}
\end{equation}
since $f$ satisfies $f(ge^{tE}mn)=e^{-t\nu}f(g)$, $m\in M_0, n\in N$.

\begin{lemm+} Let $f\in I(\nu, 0)$.  We have
  $$
  \pi(E) f(k) = I + II + III, $$
  where
$$I:=  \frac {\nu}2 \lgl k^{-1} E, E\rgl f(k), \, 
 II:= -\frac 12   \lgl k^{-1}E, F\rgl (i H f)(k), \, 
 $$
 and
 $$III:=  \left([E,  P_1( k^{-1}E)]f\right)(k).
  $$
  Here in $II$ and $III$  the expression 
  $(X f)(k)$ stands for the right differentiation 
  of $X\in \fk$ on $f$ defined in (\ref{r-diff}).
  \end{lemm+}
  \begin{proof}
    Let $X(k):= -k^{-1}E\in \fp$, $k\in K$.
    For any $K$-finite element $f\in I(\nu)$,
        \begin{equation*}
      \begin{split}
  \pi(E) f(k) &=
  \frac{d}{dt}f(\exp(-tE)k)|_{t=0}
  =  \frac{d}{dt}f(k k^{-1}\exp(-tE)k)|_{t=0}
  \\
  &=  \frac{d}{dt}f(k \exp(-t k^{-1}E))|_{t=0}
  = (X(k) f)(k)
\end{split}  \end{equation*}
is a right differentiation of the vector field  $X(k)$
on $f\in I(\nu)\subset C^\infty(G)$ evaluated at $k\in K$.
Using (\ref{Xf}) and the definition of $I(\nu)$ we find
\begin{equation}
  \label{pi-E-1}
\begin{split}  
  &\quad  \pi(E) f(k) = (X(k)f)(k)\\
  &=\frac 12 \nu   \lgl k^{-1} E, E\rgl f(k) 
-   \frac 12    \lgl k^{-1}E, F\rgl  (i H f)(k) -
\(  [ P_1( k^{-1}E), E] \)f(k)\\
& =\frac{\nu} 2 
\lgl k^{-1} E, E\rgl f(k)
-   \frac 12    \lgl k^{-1}E, F\rgl  (i H f)(k)
+\( [E,  P_1( k^{-1}E)] f\)(k)\\
&
=: I+ II +III.
\end{split}
\end{equation}
\end{proof}

In the rest of the proof we fix $f:=
\phi_{\mu; l, p}$.

\subsubsection{Expansion of $I$}

The term $I$ can be found using Proposition 
\ref{prod-exp-phi}, 
$$
I 
=
\sum_{\sigma\in \mathbb Z_2^4, \delta\in \mathbb Z_2}
I_{\sigma;  \delta}
{\phi_{\mu+\sigma; l+ \delta_1, p+ \delta_2}}
$$
with the coefficients 
    \begin{equation}
    \label{exp-fg-I}
    I_{\sigma;  \delta_1, \delta_2}: 
    =
    \frac{\nu} 2 C_{\sigma;  \delta_1, \delta_2}. 
  \end{equation}

 \subsubsection{Expansion of 
   of  $II$}
It follows from 
  (\ref{eq:1-4-sph-}), 
   (\ref{eq:H-on-f-0}) and (\ref{eq:H-on-psi}) 
 that 
\begin{equation}
  \begin{split}
    II: &= -  \frac 12 \langle k^{-1}E, F\rangle (i H f)(k)\\
    &=p 
\phi_{1^4; 1, 1}^{-}(k) 
\phi_{\mu; l, p}(k) 
+p \phi_{1^4; -1, 1}^{-}(k) 
\phi_{\mu; l, p}(k).
\end{split}
\end{equation}
We use   (\ref{eq:1-4-sph-})
and Proposition 
 \ref{prod-exp-phi} again to 
find that  the coefficients  of  $II$
  in 
 $ \psi_{\mu+\sigma; l+ \delta_1, p+\delta_2}$   are 
$$
II_{\sigma; \pm 1, 1}:= pC_{\sigma;  1, 1}, \, 
II_{\sigma; \pm 1, -1}:= -pC_{\sigma;  1, 1}, 
$$
Note that the choice of   $\pm p$ in
$II_{\sigma; \pm 1, \pm 1}$ is consistent:
$II_{\sigma; \pm 1, 1}$ has   $p$
and $II_{\sigma; \pm 1, -1}$ has   $-p$.

    \subsubsection{Expansion of 
 of  $III$}
This requires some more computations 
involving Casimir elements in $\fk$.

Let $\{u_j\}$ be an ONB 
of $\fp_1\subset (\fp, \langle \cdot, \cdot\rangle)$
as in 
Section \ref{cas-sec} with 
$[E, u_j]=z_j$ an ONB
of $(\fk_1, -\langle \cdot, \cdot\rangle)$. 
The $\fp_1$-component $P_1(k^{-1}E)=\sum_j 
\langle k^{-1} E, u_j\rangle u_j 
$,  and 
$$
[E, P_1(k^{-1} E)]
=\sum_{j}\langle 
k^{-1} E, u_j\rangle 
[E, u_j]
=\sum_{j}
\langle 
k^{-1} E, u_j\rangle z_j 
$$
for $[E, u_j]=z_j.$
Thus 
\begin{equation*}
  \begin{split}
III &=  \( [E,  P_1( k^{-1}E] f\)(k) 
  = \sum_j \langle k^{-1}E, u_j \rangle 
  (z_j f)(k) \\
&  = \sum_j \langle k^{-1}E, u_j\rangle (z_j f)(k). 
  \end{split}
\end{equation*}
We shall compute $III$
using the Casimir elements $\Cas$
and $\Cas_1$ defined in (\ref{defi-cas}).

Firstly, we claim that each term $ (k^{-1}E, u_j)$ is  obtained 
  by the right differentiation of $z_j$
  on $\langle k^{-1}E, E\rangle =
  \psi_{1^4}(k)=\langle k^{-1} E, E\rangle
  $; indeed 
\begin{equation*}
  \begin{split}
    (z_j  \psi_{1^4})(k) 
&  =    \frac{d}{dt} \langle (k\exp(tz_j))^{-1}E, E\rangle 
|_{t=0}
 =  \frac{d}{dt} \langle (\exp(-tz_j)k ^{-1}E, E\rangle 
  |_{t=0}\\
&  = \frac{d}{dt} \langle k^{-1}E, 
\exp(tz_j)  E\rangle 
|_{t=0}
=  \langle k^{-1}E, 
[z_j,  E]\rangle 
\\
&= -\langle k^{-1}E, 
[E, z_j]\rangle 
=-\langle k^{-1}E, 
u_j\rangle. 
  \end{split}
\end{equation*}
 Thus 
 $$
III=  \sum_j \langle k^{-1}E, u_j\rangle 
  (z_j f)(k) 
  = -\sum_j   (z_j  \psi_{1^4}  )  (k)   (z_j f)(k). 
  $$
  The sum $2III$
  is obtained from $\Cas_1$ by the Leibniz' rule, namely 
    \begin{equation*}
      2III 
      = - 2\sum_j (z_j \psi_{1^4}
      )( z_j f) 
  = -\Cas_1(\psi_{1^4} f) + f \Cas_1 (\psi_{1^4}) + \psi_{1^4} \Cas_1 (f). 
  \end{equation*}
Furthermore $\Cas_1 F=\Cas F-\Cas_2F-\Cas_0F =\Cas F-\Cas_2 F$
 for any $F\in L^2(K/L)$ as $\Cas_0 F=0$ by the right-invariance  of 
  $L$, so that 
  \begin{equation}
    \label{ABC-1}
    \begin{split}
&\quad III\\
&=- \frac 12 (\Cas -\Cas_2) (\psi_{1^4} f ) 
+\frac 12 f(\Cas -\Cas_2)\psi_{1^4} 
+\frac 12 \psi_{1^4} (\Cas  -\Cas_2) f.
\end{split} 
  \end{equation}

  Next all terms above can be obtained by the expansion of the
  product $\psi_{1^4}f 
  =( \phi_{1^4; 1, 1}^+ +
   \psi_{1^4; -1, 1}^+) f 
  $ as  in  (\ref{exp-fg-1})
  as sum of 
  $\phi_{\mu+\sigma; l\pm 1, p\pm 1 }:=
  \phi_{\mu+\sigma; \pm 1, \pm 1 }
  $ (dropping the index $(l, p)$ for notational convenience) 
  and each summand is 
  an eigenfunction of $\Cas -\Cas_2$, 
  $$(\Cas -\Cas_2) 
  \phi_{\mu+\sigma;  \delta }
  = (\Cas(\mu+\sigma) - \Cas_2(p +\delta)) 
  \phi_{\mu+\sigma; \delta }, 
  $$
  with
  $$\,  \Cas(\mu +\sigma) 
  =-\frac 12 \sum_{j=1}^4 (\mu_j +\sigma_j +2\rho_j) (\mu_j
  +\sigma_j),\,
  \Cas_2(p+\delta) =-2(p+\delta_2)^2.
$$
We find then
\begin{equation*}
  \begin{split}
III &=
-\frac 12 
\sum_{\sigma; \delta}
\[
\Cas(\mu +\sigma) -\Cas(\mu) -\Cas (1^4) 
  -\Cas_2(p+\delta_2)+
 \Cas_2(p) 
 +\Cas_2(1)) 
 \]
C_{\sig; \delta} \phi_{\mu+\sigma; \delta}
\\
&=\frac 12 
\sum_{\sigma; \delta}\[\sum_j(\mu_j +\rho_j)\sigma_j - \sum_j \rho_j - 4p\delta_2\]
C_{\sigma; \delta}\phi_{\mu+\sigma; \delta}
\end{split}
  \end{equation*}
 and the coefficients $III_{\sig, \delta}$ of 
$C_{\sigma; \delta}\phi_{\mu+\sigma; \delta}$ is 
$$
III_{\sig, \delta}=\frac 12 \(
\sum_j(\mu_j +\rho_j)\sigma_j - \sum_j \rho_j - 4p\delta_2\).
$$
Here we have used
the eigenvalue formula
(\ref{cas-mu})-(\ref{cas-mu-2}) 
and that
\begin{equation*}
  \begin{split}
&\quad \frac 12\(
\sum_j(\mu_j +\sigma_j +2\rho_j)(\mu_j +\sigma_j) 
-\sum_{j=1} (\mu_j +2\rho_j)\mu_j 
    -\sum_{j=1} (1 +2\rho_j)
    \) \\
    &=\frac 12 \(\sum_j(\mu_j +2\rho_j)\sigma_j +
      \sum_j \sigma_j (\mu_j +\sigma_j) 
      -\sum_{j=1} (1 +2\rho_j)
      \) \\
    &=\sum_j(\mu_j +\rho_j)\sigma_j - \sum_j \rho_j,
  \end{split}
\end{equation*}
since $\sigma_J=\pm 1$.

\subsubsection{Proof of Theorem \ref{pi-E-W}}
\begin{proof}
  Altogether
  \begin{equation}
    \label{eq:pi-nu-0-phi}
  \pi_{\nu, 0}(E)\phi_{\mu+\sigma; l, p}
=\sum_{\sigma\in \mathbb Z_2^4, \delta\in \mathbb Z_2^2}
A_{\nu; \mu; l, p}
(\mu+\sigma;  l+\delta_1, p+\delta_2) 
{\phi_{\mu+\sigma; l+\delta_1, p+\delta_2}}    
  \end{equation}
with 
$A_{\nu; \mu; l, p}
(\mu+\sigma;  l+\delta_1, p+\delta_2) $ as stated.
Here the contribution of $\nu$
is due the term $I$, 
$\sum_j(\mu_j +\rho_j)\sigma_j - \sum_j \rho_j $
is due to $II$
and $-2p \delta_2 $ is 
due to $III$ and $II$. 

The coefficients
$  C_{\nu; \mu; l, p}(\mu+\sigma;  l\pm 1, p\pm 1) $
depends only on $|l|, |p|,
|l\pm 1|, |p\pm 1|$, namely
the $2^4$ coefficients
$C_{\nu; \mu; \pm l, \pm }(\mu+\sigma;  l\pm 1, p\pm 1) $
become only $4$ possible distinct ones, and that
$ \nu +\sum_{j=1}^4 
  (\mu_j +\rho_j)\sigma_j -
  \sum_{j=1}^4\rho_j  -2 p\delta_2 $
  depends only of $p \delta_2 $. 
We may then take the even
and odd parts as in (\ref{I-pm})
and 
find the expansions
of
$\pi_{\nu}(E)\phi_{\mu+\sigma; l, p}^+$
and $\pi_{\nu, -}(E)\phi_{\mu+\sigma; l, p}^-$.
This completes the proof.
\end{proof}

\section{
  Complementary series and 
  composition series 
}
\label{comp-sect}

\subsection{Complementary series}

\begin{theo+} \label{compl-ser}
Let $\nu=\rho_{\fg} + \widetilde\nu 
=1+\frac 12 \dim{\fn_1}
+ \widetilde\nu \in \mathbb R$. 
The range of the complementary series for $I(\nu)$
    is  $|\widetilde \nu| <\widetilde \nu_0$, where 
   $$
\widetilde\nu_0=
\begin{cases}
1, \,  &\fg= \fe_{6(2)},\fe_{7(-5)}, \fe_{7(7)}, \fe_{8(-24)}, \fe_{8(8)}, 
\\
  \text{none},  &\fg= \ff_{4(4)},
                  \fe_{6(6)},  \\
\text{none},   &\fg= \fso(4, d), d \, \text{odd}\\
1, \,   &\fg= \fso(4, d), d \, \text{even}
\end{cases}.
$$
There exists no complementary series
for $I(\nu, -)$ for real  $\nu$.
\end{theo+} 
\begin{proof} 
By general theory (see e.g.~\cite{HT}) 
combined with our Theorem \ref{pi-E-W} the range is 
 determined by those $\nu\in \mathbb R$
such that $I(\nu)$ is irreducible and that
 the two linear functions
 \begin{equation}
   \label{lin-fun}
\nu +\sum_{j=1}^4 (\mu_j+\rho_j)\sigma_j -\sum_{j=1}^4 \rho_j
 -2 p\delta_2, \,
-\left(
\nu -\sum_{j=1}^4 (\mu_j+\sigma_j+\rho_j)\sigma_j  -\sum_j \rho_j
+2
(p+\delta_2)
\delta_2
\right)
\end{equation}
have the same sign. 
Writing
$\nu=\rho_{\fg} + \tilde\nu =1+\sum_j\rho_j+ \tilde\nu$ and
using Lemma \ref{rho-rho}
we see that this is precisely
the range $|\tilde\nu| <\tilde\nu_0$
where
\begin{equation*}
\begin{split}
\tilde\nu_0
&=\min_{\mu, p, \sigma, \delta_2} |1+\sum_{j=1}^4(\mu_j+\rho_j) \sigma_j -2p \delta_2  |
\\
&=\min_{\mu, p, \sigma, \delta_2} 
|1+(\mu_1+1) \sigma_1 
+(\mu_2+1+2a) \sigma_2  +
(\mu_3+1+a) \sigma_3+ 
(\mu_4+1) \sigma_4
 -2p \delta_2  |,
\end{split}  
\end{equation*}
the minimum being taken among all possible $\mu, \sigma, \delta_2$ as in the statement
of Theorem
   \ref{pi-E-W}.  If further
   $\fg\neq \ff_{4(4)}, \fso(4, d),  \fe_{6(6)}$. 
 the root multiplicity  $a$
is even and we then find $\tilde\nu_0=1$
by taking specific choices of 
 $\mu, \sigma, \delta_2$. 
 If $\fg= \ff_{4(4)},\fe_{6(6)},
 $ then $a=1$ and  we find $\tilde\nu_0=0$.
If  $\fg=\fso(4, d)$, then
 $\nu_0=0$ or $1$ depending on $d$ is odd or even, respectively.

The odd spherical representation
$I(\nu, -)$ has no complementary series since the $(K, L_1)$-types
are of the form $(\mu, p)$ with $|p|\ge 1$ and the
linear functions 
   (\ref{lin-fun}) take different values.
\end{proof}

\subsection{Reduction points and  finite dimensional 
subrepresentations}
The reduction points
are direct consequence of Theorem    \ref{pi-E-W} above.
The composition series is rather complicated 
due to the higher multiplicities.

\begin{theo+} \label{redu-pt}
Let $\nu=\rho_{\fg} + \tilde\nu 
=1+\frac 12 \dim{\fn_1}
+ \tilde\nu \in \mathbb R$.
\begin{enumerate}
\item
  The reduction 
points of $I(\nu)$ and $I(\nu, -)$ appear at the following
integer points of $\tilde\nu$, 
$$
\tilde \nu: 
\begin{cases}
\text{odd}, \,  &\fg=
\fe_{6(2)}, 
\fe_{7(-5)}, \fe_{7(7)}, \fe_{8(-24)}, \fe_{8(8)}, \\
  \text{even},  &\fg= \ff_{4(4)},
  \fe_{6(6)} \\
\text{even},   &\fg= \fso(4, d), d \, \text{odd}\\
\text{odd},   &\fg= \fso(4, d), d \, \text{even}.
\end{cases}
$$
\item 
$I(\nu)$ has a finite dimensional 
    submodule if  $\nu=-4m$, $m=0, 1, 2\cdots$, which
    is equivalent to the irreducible representation  in  $S^{2m}(\fg^{\mathbb C})$
with the highest weight vector
defined by the matrix coefficient of $(F-iH)^{2m}$,
$
f(g)=(\Ad(g^{-1}) (F-iH)^{2m}, (F-iH)^{2m}), g\in G.
$
\item 
$I(\nu, -)$  has a finite dimensional 
    submodule if  $\nu=-2m\le -2$ is an  even  integer. 
    It is equivalent to the irreducible representation  
in  $S^{m}(\fg^{\mathbb C})$
generated by the function $f(g)=B(\Ad(g) (F-iH)^{m}, (F-iH)^{m})$,
where $B$ is the Killing form extended
to   $S^{n}(\fg)$.
\end{enumerate}
\end{theo+}

\section{Minimal representation 
and  conformally 
 invariant  operators for quaternionic $\fg$
}
\label{quat-ker}

 In \cite{GW} Gross and Wallach proved
 that the minimal quaternionic
 representation $\pi^{min}$
 of  
 $$
 \fg=  \fe_{6(2)}, \fe_{7(-5)}, \fe_{8(-24)}, \fso(4, 4)$$
 appears in the induced representation
 $\pi_{\nu}$
 for   $\nu=a+2$ where $a$
 in the root multiplicity
 in  (\ref{rho-k})
 with $a=0$ for
$\fso(4, 4)$. This is done
by using the algebraic description 
of the quaternionic minimal representation 
$\pi^{min}$ and by construction
of an intertwining operator from
$\pi^{min}$
into the kernel
 of a conformal invariant differential
 operator $\Omega$ in the
 induced representation  $I({\nu})$.
The $K$-types in the minimal representation, written 
using highest weights, are 
$$
\pi^{\min}
=\sum_{n=0}^\infty 
W_{(a+2 +n, n1^3)}, \,
{(a+2 +n, n1^3)} = \frac 12((a+n)\beta_1 + n(\beta_1 + \beta_3 
+\beta_4).
$$
Here we shall give a different
proof of the result that
there exists a submodule
in $I(\nu)$
with the given $K$-types above by
using  Theorem
\ref{pi-E-W} and by the invariance
of $\Ker \Omega$.
We prove also that $\pi^{min}$ is
proper subrepresentation of the
the kernel. We shall be very brief.

Let $\{X_{\pm 2}\}$ be the $\pm 2$-root vectors of $E$ in
  (\ref{n12}).
For any $z\in \fm$ let
$\{u_j\}$ be a basis of $\fn_{1}$
and $\{u_j'\}$ a dual basis of  $\fn_{-1}$. The system
of invariant differential operators is defined as
$$
\Omega(z)
=\sum_j u_j' [X_{-2}, [z, u_j]]\in S^2(\fn_{-1}) \subset U(\fg), \,
z\in \fm
$$
as an element in the universal enveloping algebra of $U(\fg)$.
We compute $\Omega(z)$ modulo $U(\fg)\fn$
represented as a sum $f(z)E$ and
and an element in $S^2(\fk)$. 
The group $M$ acts on $\fn_{-1}$ as
a real symplectic transformation so the linear function
$\Omega(z)$ is vanishing modulo $U(\fg)\fn$
except possibly for $z\in \fl$
is in the center of $\fl$. We fix $z_0=\frac 14 (3\beta_1^\vee
-\beta_2^\vee-\beta_2^\vee-\beta_4^\vee)$.
Then $z_0$ is in the center as can be checked
by routine computation using
the root systems of $\fg$ and  $\fm$.
We recall the invariance
of $\Omega:=\Omega(z_0)$
and its component in $U(\fg^{\mathbb C})$
 modulo $U(\fg) (\fm+\fn)$.

We let   $\fg=  \fe_{6(2)}, \fe_{7(-5)}, \fe_{8(-24)}, \fso(4, d)$, $d\ge 4 
   $  even.

 \begin{lemm+}
   \label{inv-conf}
   \begin{enumerate}
\item (\cite[Theorem 5.2]{B-K-Z},
  \cite[Theorem 3.2.3, Proposition 4.7.6]{F})
   The space  $ \Ker  \Omega_{I(\nu)}$
   is    subrepresentation of $I(\nu)$.
  \item   (\cite[pp. 117-118]{GW})
Let $z_0\in \fl$
  be the center element in $\fl$ defining 
the complex structure of $(\fm, \fl)$. 
We have 
$$
\Omega
=\frac{(2d-8)}2   {E}  +
3 \(
2\Cas_{\fsu(2)} -
2\Cas_{\fu(1)}\) 
- 
\(
2\Cas_{\fm_c} - 2\Cas_{\fl_c}. 
\) 
$$
modulo $U(\fg) (\fm+\fn)$. 
  \end{enumerate}
  \end{lemm+}

\begin{prop+} \label{min-as-sub}
  Let $\fg=  \fe_{6(2)}, \fe_{7(-5)}, \fe_{8(-24)}, \fso(4, d)$, $d\ge
  4$ even and    $\nu=a+2$.
  Then
  the Harish-Chandra module
  $ I(\nu)$ has a composition
  series $\pi^{min}\subsetneq \Ker \Omega\subsetneq I(\nu)$
  with $\pi^{min}$ having  one-dimensional $K$-types $\mu=(a+n, n1^3), n\ge 0. 
$
\end{prop+}
We shall need  some more precise details
for the recurrence formula.
The  representation
$W_\mu$
of $\fk=\fsu(2)+\fm_c$ with highest weight of the form
$\mu=(\mu_1, n1^3)$
is well-understood, namely
it is the tensor product $\odot^{\mu_1}\mathbb C^2
\otimes V_{n1^3}$ with
$n1^3$ being of a central character
for the symmetric space $M_c/L_c$.
The spherical polynomials $\phi_{\mu, l}$
in this case are  the Jordan determinant \cite{Z-tams}.
In particular it implies that the highest weights appearing
in the recurrence formula $\phi_{1^4, (\pm 1, \pm 1)}\phi_{\mu, (l,
  p)}$
are  of the form $\mu=(\mu_1\pm 1, (n\pm 1)1^3).$

\begin{proof}
When acting on the space  $W_{\mu, p}
\subset  L^2(K/L_1; \chi^p) 
\subset  L^2(K/L_0)=I(\nu, 0)$
the operator 
$2\Cas_{\fsu(2)} -
2\Cas_{\fu(1)}
$ has eigenvalue 
$$
-\frac  12 (\mu_1 +2\rho_1)\mu_1 (\beta_1, \beta_1) 
+  (d\chi^p(H))^2 
= - (\mu_1 +2\rho_1)\mu_1  + p^2 
$$
and 
$
2\Cas_{\fm_c} - 2\Cas_{\fl_c}$ has eigenvalue 
$$
2\Cas_{\fm_c} - 2\Cas_{\fl_c}
= - (\mu_2 +2\rho_2)\mu_2 
- (\mu_3 +2\rho_3)\mu_3 
- (\mu_4 +2\rho_4)\mu_4 
+ p^2 
$$
since the character $d\chi^p$ on $\fl_c$ is supported on 
the center element $\frac 12 (H_{\beta_2} +
H_{\beta_3} +
H_{\beta_4})$. Thus 
$\Omega 
$ has eigenvalue 
$$
\Omega|_{W_{\mu, p}}=
(d-4)\nu 
- 3(\mu_1 + 2)\mu_1 
+ (\mu_2 +2\rho_1) \mu_2 
+  (\mu_3 +2\rho_3) \mu_3 
+ (\mu_4 +2\rho_4) \mu_4. 
$$
     Now the equation
     $   \Omega|_{W_{\mu, p}}=0$
is equivalent to
\begin{equation}\label{ome=0}
  3a (a+2) -
 3(\mu_1 + 2\rho_1)\mu_1 
+ (\mu_2 +2\rho_2) \mu_2 
+  (\mu_3 +2\rho_3) \mu_3 
+ (\mu_4 +2\rho_4) \mu_4=0, \, \mu_j= p \pmod 2 \forall j,
\end{equation}
with
$\mu=(a+n, n1^3)
$
being a special solution.
We prove now that the cyclic subrepresentation
generated by $\mu=(a, 0, 0, 0):= (a, 0)$
has the given $K$-types, i.e.
$\pi_\nu(E)^n  \phi_{(a, 0), 0}^+$
is of  $K$-type $
(a+n, n1^3)$.
It follows from Theorem    \ref{pi-E-W}
for $\nu=a+2$ we have
$\pi_\nu(E)  \phi_{(a, 0), 0}^+$
is a linear combination of
$\phi_{(a+1, 1^3); \pm 1}^+$ and
$\phi_{(a-1, 1^3); \pm 1}^+$
with the coefficient of
$\phi_{(a-1, 1^3); \pm 1}^+$ vanishing;
in other words $\pi_\nu(E)  \phi_{(a, 0), 0}^+$
is of $K$-type $
(a+1, 1^3)$. 
Suppose this
is true for a fixed $n\ge 1$. 
The $K$-type in
$\pi_\nu(E)^{n+1}  \phi_{(a, 0), 0}^+$
is then of the form
$$
\mu=(a +n +m_1, (n +m_21)1^3), \quad m_1=\pm 1, m_2=\pm 1 
$$
for $\fg\neq \fso(4, d)$,
$$
\mu=(a  +n+m_1, n +m_2, (n+m_3)1^2), \quad  m_1=\pm 1, m_2=\pm 1, 
$$
for $\fg= \fso(4, d)$, $d>4$
and 
$$
\mu=(a +n +m_1, n +m_21, n+m_3, n+m_4),  m_1=\pm 1, m_2=\pm 1,
m_3=\pm 1, m_4=\pm 1, 
$$
for $\fg= \fso(4, 4)$. Inserting these $\mu$ into
the equation
(\ref{ome=0}) we find that $\mu=(a+n+1, (n+1)1^3)$
or $\mu=(a+n-1, (n-1)1^3)$ are the only possible
solutions among these possible $\mu$.
This proves that there is a subrepresentation
$\pi^{min}\subseteq \Ker \Omega$.

The equation (\ref{ome=0})
is independent of $p$ so that
for each solution $(a+n, n, n, n)$
there is a multiplicity of $[\frac n 2]$
counting those $p$, $ p=n \pmod 2$.
However $\pi^{min}$
is of multiplicity one. This proves
that $\pi^{min}\subsetneq \Ker \Omega$, 
the proper inclusion
$\Ker \Omega\subsetneq  I(\nu)$
is trivial. This finishes the proof.
\end{proof}

\begin{rema+}     Note here that we have
  used Lemma \ref{inv-conf} $\Ker \Omega$ is
  a subrepresentation of $I(\nu)$.
 It might be possible to prove that
 $I(\nu)$ has a subrepresentation with one-parameter
 $K$-type generated by the $K$-type $(a, 0, 0, 0)$
 by finding a closed
combinatorial formula for the functions $\pi_\nu(E)^n\phi_{(a, 0, 0, 0), 0}^+$.
   For the Lie algebra  $\fg_2$
   Kable \cite{Ka} has found
   such a formula. Also
the  equation (\ref{ome=0})
is of the form $Q(\mu)=c$
where $Q$ is a quadratic form
with integral coefficients of signature $(3, 1)$
and it has solution except the line
$(a+n, n, n, n)$ by the Borel
density theorem \cite{WM}.
\end{rema+}

\appendix

\section{Spherical polynomials of type $(l, p)$ and
tables of $(\fg, \fk)$ and $(\fm, \fl)$}
\label{appen-sph-lp}
The theory of spherical functions
 on  general compact
 symmetric spaces is well developed; see e. g.
 \cite{He2}. For 
 a Hermitian
 symmetric space
 there are
also  spherical functions corresponding
to holomorphic line bundles
see \cite{Sch, Shi}.
To avoid excessive
 notation we keep our 
setup and consider only
the compact Hermitian symmetric space
 $K/L_1=\mathbb P^1\times M_c/L_c$ 
above.
These spherical functions
$\phi$  \cite{Sch}
as functions on  $K$
 are defined as eigenfunctions of
 invariant differential operators
 such that
 $\chi^l(k_1)\chi^l(k_2)\phi
 (k_1kk_2)= \phi(k) 
 , k\in K,  k_1, k_2\in L_1$
for the  character $\chi^l$. However
for our study of induced representations
 we need product formulas for spherical functions 
 of type $(l, p)$, i.e.
 eigenfunctions
 $\phi$ such that $
 \chi^l(k_1) \chi^p(k_2)\phi(k_1 gk_2)=
 \phi(g)$, $k_1, k_2\in L_1$, 
with $l$ and $p$  satisfying  certain compatible conditions.
We shall prove that
the product formulas for spherical functions 
of type $(l, p)$ can be obtained
from those of diagonal type $(l, l)$.

\subsection{Spherical polynomials of type $(l, p)$
on  Hermitian symmetric space $K/L_1$}

Let $p=l \pmod 2$ be integers.  We call \cite{Sch, Shi} a function 
$f\in C^\infty(K)$ of type $(l, p)$ under $L_1$ if
\begin{equation}
  \label{eq:p-l-type}
  \chi^l(k_1)
  \chi^p(k_1)f(k_1kk_2) 
=  f(k), \, k\in K, k_1, k_2\in L_1.
\end{equation}
To clarify the construction of
spherical polynomials $\phi_{\mu; l, p}$
of type $(l, p)$ we consider
a unitary representation
$(W_{\mu}, \tau_\mu, K)$
 with highest weight 
 $\mu=\frac 12 (\mu_1\beta_1+
 \mu_2\beta_2+\mu_3\beta_3 +\mu_4\beta_4)$.
 with $W_{\mu, p}$ being
  realized as in
 $W_{\mu, p}\subset L^2(K/L_1, \chi^p)$ by
 (\ref{L2KL+-0}).
 In particular for any integer $l$
 satisfying the same condition as $p$ 
there exists a unique unit vector up to scalars
 $v_l=v^{\mu}_l\in W_{\mu}$ such that
 $\tau(k)v_l=\chi^l(k) v_l, k\in L_1$.

\begin{defi+}
We define the spherical polynomial $\phi_{\mu; l, p}$
of type $(l, p)$ by
$$
\phi_{\mu; l, p} (k)
=\langle
\tau_\mu(k^{-1})v_{l}, v_{p}\rangle, \, k\in K.
$$
\end{defi+}
Then  $\phi_{\mu; l, p} $ satisfies
the transformation rule
$$
\chi(k_1)^l\chi(k_2)^p 
\phi_{\mu; l, p} (k_1kk_2)=
\phi_{\mu; l, p} (k), \, k\in K, k_1, k_2\in L_1.
$$
We note that there is an ambiguity
in the notation
$\phi_{\mu; l, p} (g)$
as it depends also on the choice
$(v_{l}, v_{p})$, whereas 
$\phi_{\mu; p} (g):=\phi_{\mu; p, p} (g)
$ is independent of the choice of
a unit vector $v_p$.  We fix
 a choice of $\{v_p\}$
for all $p$ 
in the fixed representation space $W_\mu$ such that
$ \tau_{\mu}(w_0) v_p = v_{-p}$. With this
choice we have
\begin{equation}
  \label{eq:conv-lp}
  \phi_{\mu; l, p}^{\pm} \in W^{\pm}_{\mu; p},
  \phi^{\pm}_{\mu; l, p}(k):=
  \phi_{\mu; l, p}(k) \pm    \phi_{\mu; l, p}(k w_0), 
  k\in K, \, p\ge 0
\end{equation}
with the interpretation that
  $\phi^{-}_{\mu; l, 0}=0$.

\subsection{Recursion formula for the
  product  $\phi_{\mu; l}  \phi_{1^4; 1}
  $
of spherical polynomials of type $(l, l)$}

The Harish-Chandra theory
for spherical functions
can also be extended to
spherical functions of type $(l, l)$;
see \cite{Sch}. 
The symmetric space $K/L_1$  is a product
so it the $c$-function
\begin{equation}
  \label{eq:c-funct}
  c(-i(\mu+\rho); l)=
  c^{SU(2)}(-i(\mu_1+\rho_1); l) 
  c^{M_c}(-i(\mu_2+\rho_2),
  -i(\mu_3+\rho_3),
  -i(\mu_4+\rho_4),  l).
\end{equation}

Recall that $\fp^{\mathbb C}
$
has highest weight
$\gamma=\frac 12(\beta_1+\beta_2+\beta_3+\beta_4):=1^4$.
Now  $\Ad(w_0)E_\gamma= E_{-\gamma}$
and this coincides with our construction
of $\phi_{1^4, \pm 1, \pm 1}$ above, so that
 $\phi_{1^4, \pm 1, \pm 1}
=\langle k^{-1} E_{\pm \gamma}, E_{\pm \gamma}\rangle$.
We shall need an expansion for the spherical 
$\phi_{1^4; \pm, \pm} \phi_{\mu; l, p}$
and we consider the case
$\phi_{1^4; 1} \phi_{\mu; l}$
first by using the known method \cite{Vretare, Z22}.

By exactly the same arguments as in 
\cite{Vretare} for $p=0$ and
\cite[Lemma 3.5]{Z22} for general $p$
that 
\begin{equation}
  \label{eq:prod-ll}
  \phi_{1^4, \pm 1}\phi_{\mu;  p}
  =\sum_{\sigma\in \mathbb Z_2^4}  C_{\mu, l}(\mu +\sigma;  p\pm 1) 
  \phi_{\mu +\sigma;  p\pm 1}
\end{equation}
where 
$$
C_{\mu, p}(\mu +\sigma;  p\pm 1) 
=\frac{
  c(-i(1^4+\rho), 1)c(-i\sigma(\mu+\rho), p ) 
}
{
 c(-i(\sigma(\mu) + 1^4+ \rho), p\pm 1 ) 
}.
$$
Here with some abuse 
of notation we have denoted 
$\sigma\in \mathbb Z_2^4=\{\pm\}^4$
both as Weyl group elements 
and as a weight $\sigma=(\sigma_1, 
\sigma_2, 
\sigma_3, 
\sigma_4)=\frac 12(\sigma_1\beta_1 +
\sigma_2\beta_2 +
\sigma_3\beta_3 
\sigma_4\beta_4). 
$
Moreover $C_{\mu, p}(\mu +\sigma;  p\pm 1) >0$ as 
far as $(\mu; p)$ and $(\mu+\sigma; p\pm 1)$ satisfy 
the dominate conditions 
(\ref{rank4}), 
 (\ref{so-4-d}), (\ref{so-4-4}).

\subsection{Recursion formula for
  $\phi_{1^4; \pm 1, \pm 1}\phi_{\mu; l, p}$
  for $\phi_{\mu; l, p}$ of general type $(l, p)$
}
 \begin{prop+}
 \label{prod-exp-phi}
There exist positive constants $C_{\pm 1, \pm 1}(\mu +\sigma; l\pm 1, p\pm 1)$
such that
  \begin{equation}
    \label{exp-fg-1}
{\phi_{1^4; \pm 1, \pm 1}}   
 \phi_{\mu; l, p}
=    \sum_{\sigma\in \mathbb Z_2^4} C_{\mu, p}(\mu +\sigma; l\pm 1, p\pm 1)
\phi_{\mu +\sigma;  l\pm 1, p\pm 1}
+    \sum_{\sigma} C_{\mu, p}(\mu +\sigma; l\pm 1, p\pm 1)
\phi_{\mu+\sigma; l\pm 1, p\pm 1}.
\end{equation}
where $\pm 1$ are taken independently
in the LHS and $\pm 1$ are chosen accordingly
in the RHS. Moreover the constants
depend only on $
|l|, |p|,
|l\pm 1|, |p\pm 1|$.
\end{prop+}

\begin{proof} 
  The expansion   (\ref{eq:prod-ll}) for $\phi_{1^4; 1}$
  is equivalent to
\begin{equation}
  \label{l-l-case}
  \langle \tau_{1^4}\otimes
  \tau_{\mu} (k^{-1}) (
    E_\gamma\otimes v_{p}, 
    E_\gamma\otimes v_{p}\rangle
=\sum_{\sigma\in \mathbb Z_2^4} C_p(\sigma)
 \langle \tau_{\mu +\sigma} (k^{-1}) v_{p+1}^{\mu +\sigma}, v_{p+1}^{\mu +\sigma}\rangle
\end{equation}
for some positive constants $C(\sigma)$.
We rephrase this as tensor product
decomposition of
$\fp^+\otimes W_{\mu}$.
So let
$$\fp^+\otimes W_{\mu}
\sim \sum_{\sigma} m_\sigma W_{\mu +\sigma} + \text{Rest}
$$
be the irreducible decomposition
where $m_\sigma$ is the multiplicity
and the  term Rest
is a  sum of irreducible components containing
non one-dimensional $L_1$-type.
(One can prove that $m_\sigma=1$ but
we shall not need it here.)
This is  equivalent to 
$$
E_\gamma\otimes v_{p}
=\sum_{\sigma} a_p(\sigma)
v_{p+1}^{\mu+\sigma} + \text{rest}
$$
for some
constants $a(\sigma)$ and {\it some unit vectors} $
v_{p+1}^{\mu+\sigma} \in W_{\mu +\sigma}
$  with highest weight  $ {\mu+\sigma}$.
Comparing this with (\ref{l-l-case})
and use Schur orthogonality we
see that
$|a_p(\sigma)|^2
=C_p(\sigma)$. We can thus
take  $a_p(\sigma)=
\sqrt{C_p(\sigma)}$
by a choice of   $v_{1+l}^{\mu+\sigma}$.
Computing matrix coefficients
of $E_\gamma\otimes v_l$ with
$E_\gamma\otimes v_p$ 
we find 
  \begin{equation}
    \label{l-p-case}
      \begin{split}
\phi_{1^4; 1, 1}(k) \phi_{\mu; l, p}(k) 
&=\langle \tau_{1^4}(k^{-1})\otimes \tau_{\mu}(k^{-1})
(E_\gamma
\otimes v_{l}^{\mu}), 
E_\gamma \otimes v_{p}^{\mu}\rangle 
\\
&=\sum_{\sigma}
\sqrt{C_l(\sigma)} 
\sqrt{C_p(\sigma)} 
\langle 
\tau_{\mu+\sigma}( k^{-1}) 
v_{l+1}^{\mu+\sigma}, v_{p+1}^{\mu+\sigma}
\rangle\\
&=\sum_{\sigma}
\sqrt{C_l(\sigma)} 
\sqrt{C_p(\sigma)}  
{\phi_{\mu+\sigma; l+1, p+1}}(k).
\end{split}
\end{equation}
The same proof applies 
to all  products
$\phi_{1^4; \pm 1, \pm 1}(k) \phi_{\mu; l, p}(k)$.
This proves our claim.
\end{proof}

\subsection{Tables}

We give a list of the symmetric
pairs $(\fg, \fk)$ studied in
this paper with the corresponding subalgebras $\fm$
in the Heisenberg parabolic $\fm+\fa+\fn$.

\begin{table}[!h]
      \begin{center}
\begin{tabular}
{
  | l| c |c |c  |}
  \hline 
  $\fg$  &$\fk$ & $\fm$ & $\fl$ 
\\
  \hline 
    $\fso(4, d), d\ge 4$
  &
     $       \fsu(2)+ \fsu(2) +\fso(d) $
&
             $\fso(2, d-2) +\fsl(2, \mathbb R)$
                  &
$\fso(2)^2 + \fso(d-2)$
  \\
  \hline 
  $\fe_{6(2)}$
  &
    $\fsu(2) + \fsu(6)$
    &
 $\fsu(3, 3)$
          &
                            $ \fsu(3) +\fu(3)$
  \\
  \hline
  $\fe_{7(-5)}$
  &             $ \fsu(2)+
               \fspin(12) $                 
          &
            $\fso^\ast(12)$
                 &
                               $\fu(6)$
  \\
  \hline 
  $\fe_{8(-24)}$
          &
            $\fsu(2)+\fe_7$
                 &
                   $\fe_{7(-25)}$
                 & $\fe_6 +\mathbb R$
  \\
  \hline   
  $\ff_{4(4)}$                      
  &
    $\fsu(2)+\fsp(3)$
    &
                 $\fsp(3, \mathbb R)$
                 &
                   $\fu(3)$
  \\
  \hline
    \hline 
  $\fe_{6(6)}$
  &
    $\fsp(4) $
    &
 $\fsl(6, \mathbb R)$
          &
                            $ \fso(6) =\fsu(4)$
  \\
  \hline
  $\fe_{7(7)}$
  &             $ \fsu(8)$                 
          &
            $\fso(6, 6)$
                 &
                   $\fso(6)+\fso(6) =\fsu(4)+\fsu(4)
                   $
  \\
  \hline 
  $\fe_{8(8)}$
          &
            $\fso(16)$
                 &
                   $\fe_{7(7)}$
                 & $\fsu(8)$
                     \\
\hline
\end{tabular}
\vskip0.20cm 
 \caption 
 {List of non-Hermitian   Lie algebras $\fg$
   and their Heisenberg parabolic
   Lie algebras $\fa+\fm+\fn$, $\fn=\fn_2+\fn_1$
   studied in this paper; the first five
   cases of $(\fg, \fk)$ are quaternionic
   symmetric pairs of  real rank $4$,  the last
   three cases of real rank $6, 7, 8$ and are not quaternionic.}
\label{tab:1}
      \end{center}
    \end{table}
\begin{table}[!h]
      \begin{center}
\begin{tabular}
{ |l | c| c |}
  \hline 
  $\fg$
  & $a$
  & $d=\dim_{\mathbb C} K/L_1=\frac 12 \dim_{\mathbb R}\fn_1 $
\\
  \hline 
  $\fso(4, d), d\ge 4$
&  $d-4$ & $d$
\\
  \hline 
  $\fe_{6(2)}$
      & 2 
                                 &          10 
  \\
  \hline 
              $\fe_{7(-5)}$
        &
          4
                          &16
  \\
  \hline 
            $\fe_{8(-24)}$   &     8
       &28
  \\
  \hline 
  $\ff_{4(4)}$                      
                           &1
                           &7                    
  \\
  \hline
    \hline
    $\fe_{6(6)}$
      & 1 
                                 &          10 
  \\
  \hline 
              $\fe_{7(7)}$
        &
          2
                          &16
  \\
  \hline 
            $\fe_{8(8)}$   &     4
       &28
  \\
  \hline 
\end{tabular}
\vskip0.20cm 
 \caption 
 {The root multiplicity $a$
and dimension for the compact Hermitian
   symmetric pair $(K, L_1)$,
   $(\fk, \fl_1), \fl_1=\fu(1)+\fl$,
   all with $\text{rank}K/L_1=4$. The  half sum of positive
   roots $\rho_{\fg}
   = 1+ \dim_{\mathbb C} K/L_1=1+d$, $d=3a +4$
 for quaternionic  $\fg\ne \fso(4, d)$
 and $d=4+6a$ for $\fg=\fe_{6(6)},
 \fe_{7(7)},
 \fe_{8(8)}.
 $
 }
\label{tab:2}
      \end{center}
    \end{table}

\newpage

\def\cprime{$'$} \newcommand{\noopsort}[1]{} \newcommand{\printfirst}[2]{#1}
  \newcommand{\singleletter}[1]{#1} \newcommand{\switchargs}[2]{#2#1}
  \def\cprime{$'$} \def\cprime{$'$} \def\cprime{$'$} \def\cprime{$'$}
\providecommand{\bysame}{\leavevmode\hbox to3em{\hrulefill}\thinspace}
\providecommand{\MR}{\relax\ifhmode\unskip\space\fi MR }
\providecommand{\MRhref}[2]{%
  \href{http://www.ams.org/mathscinet-getitem?mr=#1}{#2}
}
\providecommand{\href}[2]{#2}

\end{document}